\documentclass[a4paper,12pt]{scrartcl}

\usepackage[utf8]{inputenc}
\usepackage[english]{babel} 
\usepackage{amsmath,bm}
\usepackage{amssymb}
\usepackage{amsthm}
\usepackage{amsfonts}
\usepackage[mathscr]{eucal}
\usepackage{algorithm}
\usepackage{algorithmic}
\usepackage{graphicx}
\usepackage{subfigure}
\usepackage{url}

\newtheorem{definition}{Definition}

\newtheorem{remark}{Remark}
\newtheorem{proposition}{Proposition}

\newtheorem{theorem}{Theorem}

\newcommand{\N}{\mathcal{N}}
\newcommand{\Y}{\mathcal{Y}}
\newcommand{\U}{\mathcal{U}}
\newcommand{\W}{\mathcal{W}}
\newcommand{\Uad}{\mathcal{U}_\mathsf{ad}}
\newcommand{\Zadeps}{\mathcal{Z}^\varepsilon_\mathsf{ad}}

\newcommand{\bx}{{\bm x}}
\newcommand{\bs}{{\bm s}}
\newcommand{\bn}{{\bm n}}
\newcommand{\ua}{u_\mathsf{a}}
\newcommand{\ub}{u_\mathsf{b}}
\newcommand{\ya}{y_\mathsf{a}}
\newcommand{\yb}{y_\mathsf{b}}
\newcommand{\wa}{w_\mathsf{a}}
\newcommand{\wb}{w_\mathsf{b}}
\newcommand{\AWon}{\mathscr A_1^\W}
\newcommand{\AWtw}{\mathscr A_2^\W}
\newcommand{\AWth}{\mathscr A_3^\W}
\newcommand{\AWfo}{\mathscr A_4^\W}
\newcommand{\AWfi}{\mathscr A_5^\W}
\newcommand{\AWsi}{\mathscr A_6^\W}
\newcommand{\AaU}{{\mathscr A_{{\mathsf a}}^\U}}
\newcommand{\AbU}{{\mathscr A_{{\mathsf b}}^\U}}
\newcommand{\AaW}{{\mathscr A_{\mathsf a}^\W}}
\newcommand{\AbW}{{\mathscr A_{\mathsf b}^\W}}
\newcommand{\AaVT}{{\mathscr A_{\mathsf a,2}^\W}}
\newcommand{\AbVT}{{\mathscr A_{\mathsf b,2}^\W}}
\newcommand{\IU}{{\mathscr I^\U}}
\newcommand{\IW}{{\mathscr I^\W}}
\newcommand{\IVT}{{\mathscr I^\W_2}}
\newcommand{\UNj}{\mathcal{U}_N^{j}}
\newcommand{\WNj}{\mathcal{W}_N^{j}}
\newcommand{\ZadepsNj}{\mathcal{Z}^{\varepsilon,j}_{N,\mathsf{ad}}}
\newcommand{\ZadepsellNj}{\mathcal{Z}^{\varepsilon,\ell,j}_{N,\mathsf{ad}}}
\newcommand{\Zadepsinfj}{\mathcal{Z}^{\varepsilon,j}_{\infty,\mathsf{ad}}}

\def\N{\mathbb{N}}

\usepackage{color}

\newcommand{\dagfootnote}[1]{
	\let\oldthefootnote=\thefootnote
	\setcounter{footnote}{1}
	\renewcommand{\thefootnote}{\fnsymbol{footnote}}
	\footnote{#1}
	\let\thefootnote=\oldthefootnote
}

\newcommand{\ddagfootnote}[1]{
	\let\oldthefootnote=\thefootnote
	\setcounter{footnote}{2}
	\renewcommand{\thefootnote}{\fnsymbol{footnote}}
	\footnote{#1}
	\let\thefootnote=\oldthefootnote
}

\allowdisplaybreaks

\DeclareOldFontCommand{\rm}{\normalfont\rmfamily}{\mathrm}
\begin{document}
\title{Parallelized POD-based Suboptimal Economic Model Predictive Control of a State-Constrained  Boussinesq approximation}

\author{Julian Andrej\footnote{Lawrence Livermore National Laboratory, 7000 East Ave., Livermore, 94550, USA.} , Lars Gr\"une\footnote{University of Bayreuth, Universit\"atsstra\ss e 30, 95447 Bayreuth, Germany.} , Luca Mechelli\footnote{University of Konstanz, Universit\"atsstra\ss e 10, 78464 Konstanz, Germany.} , Thomas Meurer\footnote{Chair of Automatic Control, Kiel University, Kaiserstra\ss e 2, 24143, Kiel, Germany.} ,\\Simon Pirkelmann$^\dag$ and Stefan Volkwein$^\ddag$}
\date{}
\maketitle
\begin{abstract}
	Motivated by an energy efficient building application, we
	want to optimize a quadratic cost functional subject to
	the Boussinesq approximation of the Navier-Stokes equations
	and to bilateral state and control constraints. Since the
	computation of such an optimal solution is numerically
	costly, we design an efficient strategy to compute a
	sub-optimal (but applicationally acceptable) solution with
	signifcantly reduced computational effort. We employ an
	economic Model Predictive Control (MPC) strategy to obtain a
	feedback control. The MPC sub-problems are based on a
	linear-quadratic optimal control problem subjected to mixed
	control and state constraints and a convection-diffusion
	equation, reduced with proper orthogonal decomposition. To
	solve each sub-problem, we apply a primal-dual active set
	strategy. The method can be fully parallelized, which enables
	the solution of large problems with real-world parameters. 
\end{abstract}

\section{Introduction}
Coupled heat-convection phenomena like those occurring in heating, cooling and air-conditioning (HVAC) of residual buildings can be accurately modeled by the incompressible Navier-Stokes equations; see \cite{AT90,Tem84,Tri77}, for instance. These equations describe how the physical quantities temperature, air velocities and pressure are connected and how they influence each other. While the underlying physical principles are completely described by the equations, working with them directly for the purpose of control of HVAC processes has turned out to be challenging due to the high computational complexity. Instead, simplified models of varying granularity are used in control applications.

One common simplification that still retains much of the relevant physical relations in the HVAC regime is the Boussinesq approximation \cite{Tem84,Tri77}. The general idea is to reduce the coupling between heat and airflow by only considering bouyancy induced changes of the fluid velocity. More specifically, we only consider density variations along the gravity field which lead to, e.g., the rise of warmer air. Although we obtain a simplified model, this simplification is not enough to prevent high computational cost for computing its solution. Many degrees of freedom (DOFs) are in fact necessary to accurately capture the involved dynamics, especially when dealing with real-world parameters; see, e.g., \cite{And19}. 

One viable approach is to apply model order reduction (MOR). This field was extensively developed during the last two decades and several techniques were proposed. Among all, in case of computational fluid dynamics, proper orthogonal decomposition (POD) seems the one to guarantee best performances for nonlinear dynamical systems; see, e.g., \cite{BMQR15,HLBR12,SR18}. In the particular case of the Boussinesq approximation, the POD-based MOR approach is also certified by an a-posteriori error analysis \cite{Rav2000,Rav2011}, that guarantees the validity of approximation for the computed reduced order model. Due to the request of computing many PDE solutions, MOR techniques are naturally applied in the field of optimal control. In \cite{HU07}, for example, a POD-based model predictive control (MPC) algorithm is considered to approximate a feedback solution to the Boussinesq equation. 

MPC is a well-established model-based control method, sometimes also called Receding Horizon Control; cf. \cite{GrPa17,RMD17}. The key idea is to decompose an infinite time horizon optimal control problem into several finite horizon sub-problems and, therefore, to optimize predictions of the model behavior in relatively short time. After a sub-problem is solved, the initial part of the computed optimal control is stored and used as feedback control, the time horizon is shifted and it is possible to update the problem parameters before the procedure is iterated. This technique is particularly useful in HVAC applications. It permits, in fact, to merge predictions, such as weather forecast, with system measurements \cite{OPJG12,SOCP11}. Let us mention that POD-based MPC is considered also in \cite{FFV19,GU14,Mec19,MV19} for different models. Furthermore, we are particularly interested in economic MPC; cf. \cite{AmRA11,ELC17,FaGM18,GP17,GP20,Pir20}. In this particular branch of MPC, more general cost functionals than those employing tracking type costs are allowed. In our specific application, the cost functional depends only on the control variables, without any dependence on the state variable. The focus is then to keep the system trajectory inside pre-defined bounds with less control effort as possible and it deals perfectly with the usual requests of HVAC applications.

This work aspires to unify the above mentioned contributions in a reliable method, in order to compute a suboptimal solution to an optimal control problem subject to the Boussinesq approximation and bilateral state and control constraints. In particular, starting from \cite{And19,Mec19,Pir20}, we propose a POD-based economic MPC scheme which considers linear-quadratic subproblems governed by convection-diffusion equations with a virtual control approach \cite{KR09,Mec19} and possible updates of the reduced order model through parallel solution of the Boussinesq approximation \cite{And19}. Furthermore a primal-dual active set strategy (PDASS) \cite{HIK02,IK03} is considered to solve each sub-problem. The goal is to have a robust and reliable method that speeds up the computational time.

In Section~\ref{sec:setting}, we present the optimal control problem. Section~\ref{sec:pdass} describes the virtual control approach and the PDASS for linear-quadratic optimal control problems, while Section~\ref{sec:POD} contains a brief explanation of the POD method and the a-posteriori error estimate. The POD-based economic MPC algorithm is explained in Section~\ref{sec:MPC} and further details about the implementation are given in Section~\ref{sec:parallel}. Finally, we test our algorithm for a simplified problem in Section~\ref{sec:num_test}.

\section{The setting of the optimal control problem}
\label{sec:setting}
The Boussinesq approximation of the Navier-Stokes equations is described by the following PDE:
\begin{subequations}
	\label{eq:Bous}
	\begin{align}
	\frac{\partial{\mathbf{v}}}{\partial t} + {\mathbf{v}}\cdot \nabla {\mathbf{v}}	&= - \frac{1}{\rho} \nabla p	+ \nu \Delta {\mathbf{v}} - \textbf{g} \alpha (	{y}	- \tilde{y}) && \text{in } Q:=(0,T)\times\Omega, \label{eq:pde-navier-stokes}\\
	\nabla \cdot {\mathbf{v}}&= 0 && \text{in } Q, \label{eq:pde-continuity}\\
	\mathbf{v}(0)& = \mathbf{v}_\circ, && \text{in } \Omega, \\
	\frac{\partial 	y}{\partial t} + {\mathbf{v}}\cdot \nabla y	&= \frac{\kappa}{\rho c_p}\,\Delta y	&& \text{in } Q,	\label{eq:pde-heat-equation} \\
	y(0) & = y_\circ && \text{in } \Omega, \label{eq:pde-heat-equation-initguess}
	\end{align}
	where $\mathbf{v}:Q\to\mathbb{R}^{2}$ is the \emph{air velocity}, $p:Q\to\mathbb{R}$ stands for the \emph{air pressure},	and $y:Q\to\mathbb{R}$ is the \emph{temperature}, $\Omega\subset\mathbb{R}^2$ denotes an open and bounded domain with sufficiently smooth boundary, $T>0$ is the time horizon, $\textbf{g}\in\mathbb{R}^2$ is the \emph{gravitational acceleration}, $\mathbf{v}_\circ$, $y_\circ$ are given initial conditions and $\tilde{y}\in\mathbb{R}$ is the reference temperature at which we measure:
	\begin{itemize} 
		\item the (constant) density $\rho>0$,
		\item the (constant) kinematic viscosity $\nu>0$, 
		\item the (constant) coefficient of thermal expansion $\alpha>0$, 
		\item the (constant) thermal conductivity $\kappa>0$,
		\item the (constant) isobaric specific heat $c_p>0$ 
	\end{itemize}
	of the studied fluid (air in our case). 
	\begin{figure}
		\centering
		\includegraphics[height=50mm]{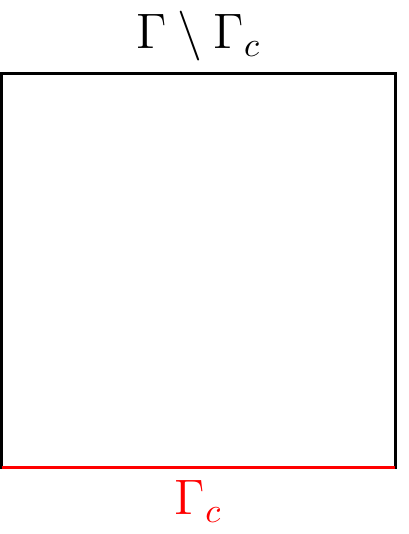}
		\caption{Domain $\Omega$ in case of underfloor heating system. \label{fig:domain}}
	\end{figure}
	In addition, we consider time-varying boundary conditions of the following form
	\begin{align}
	v & = 0 && \text{ on } (0,T)\times \Gamma, \\
	-\kappa\frac{\partial y}{\partial\bn} &= \gamma (y - y_{\text{out}}) && \text{ on }  (0,T)\times \Gamma \setminus \Gamma_\mathsf{c} ,\label{eq:pde-heat-equation-outside}
	\\
	-\kappa\frac{\partial y}{\partial\bn} &= \gamma_\mathsf{c} (y - u) && \text{ on } (0,T)\times \Gamma_\mathsf{c} ,\label{eq:pde-heat-equation-control}	 					
	\end{align}
\end{subequations}
where $\Gamma= \partial\Omega$ and $\gamma,\gamma_c>0$ are fixed parameters. The function $u\in \U= L^2(0,T)$ is a control input placed at a predefined part of the domain's boundary, namely $\Gamma_\mathsf{c}$, and it represents a device for controlling the temperature in the room, e.g. a radiator or an underfloor heating system; cf. Figure~\ref{fig:domain}. The quantity $y_{\text{out}}\in L^2(0,T)$ represents the outside temperature.

We are interested in the weak solution to \eqref{eq:Bous}. For more details of the straightforward variational formulation for \eqref{eq:Bous} we refer to \cite{AT90,Tem84} and \cite{HU07} for more details.

The goal is to minimize the control effort while keeping the temperature $y$ inside desired bounds, i.e. we are interested in minimizing
\begin{subequations}
	\label{eq:opt_cont_prob}
	\begin{equation}
	\mathcal{J}(u)= \frac{1}{2} \|u\|^2_{\U}
	\end{equation}
	subject to (s.t.) the weak formulation of \eqref{eq:Bous} and the inequality constraints
	\begin{align}
	&\ua(t) \leq u(t) \leq \ub(t) && \text{ a.e. in } (0,T) \label{control_constraints},\\
	& \ya(t,\bx) \leq y(t,\bx) \leq \yb(t,\bx) && \text{ a.e. in } Q. \label{state_constraints}
	\end{align}	
\end{subequations}
Here $y$ is the (weak) solution of \eqref{eq:Bous}, $\ua,\ub\in \U$ given control constraints and $\ya,\yb\in C(Q)$ given state constraints. It follows from \cite[Chap.~III,\,\S3]{Tem84} that \eqref{eq:Bous} admits a unique weak solution
\begin{align*}
(\mathbf{v},y)\in \Y_\mathsf B=W(0,T;\tilde V)\times W(0,T;V),
\end{align*}
where $V=H^1(\Omega)$ and the Hilbert space $\tilde V$ is the closure of the space $\{v\in C^\infty_0(\Omega)^2\,|\,\nabla\cdot v =0\}$ with respect to the norm $\|\cdot\|_{V\times V}$. Moreover, we have $W(0,T;\tilde V)=L^2(0,T;\tilde V)\cap H^1(0,T;\tilde V')$, and the set $W(0,T;V)$ is defined analogously; see \cite[Chap.~XVIII]{DL00} for more details. Thus, considering the control-to-state map 
\begin{align*}
\mathcal S_\mathsf B: \U \to \Y_\mathsf B,\quad u\mapsto \mathcal S_\mathsf B(u)=(\textbf{v}(u),y(u))\in\Y_\mathsf B
\end{align*}
and the admissible set 
\begin{align*}
\U_{\mathsf B,\mathsf{ad}}= \big\{ u\in\U\,\big|\,u\text{ and }y=\mathcal S_\mathsf B(u)\text{ satisfy \eqref{control_constraints} and \eqref{state_constraints}, respectively}\big\},
\end{align*}
problem \eqref{eq:opt_cont_prob} can be rewritten as the following purely control-constrained optimization problem
\begin{equation}
\tag{$\mathbf{\hat P}_\mathsf B$}
\label{eq:opt_cont_prob_in_u_bous}
\min \mathcal{J}(u)\quad\text{s.t.}\quad u\in\U_{\mathsf B,\mathsf{ad}}.
\end{equation}		
We call \eqref{eq:opt_cont_prob_in_u_bous} the \emph{reduced problem} because -- in contrast to \eqref{eq:opt_cont_prob} -- it is an optimization problem in the control variable only. Note that \eqref{eq:opt_cont_prob_in_u_bous} admits a global solution $\bar u_\mathsf B\in\U_{\mathsf B,\mathsf{ad}}$ provided $\U_{\mathsf B,\mathsf{ad}}$ has non-empty interior; cf., e.g., \cite{HPUU2009}. Throughout the work a bar indicates optimality.

While the Boussinesq approximation already reduces numerical effort for simulation of the fluid processes compared to the full compressible Navier-Stokes equations, the nonlinearity induced by the two-way coupling of temperature and velocity still makes it computationally infeasible, when it is used, e.g., as the prediction model in an MPC scheme. This is because each step of the MPC algorithm typically requires many evaluations of the prediction model, which would be too costly for an online application of the method. 

For that reason we elaborate a simplified strategy that aims at decreasing the computational cost significantly. We introduce several simplifications from a model and optimal control perspectives. At first, we notice that the buoyancy effect in the Boussinesq model has a big impact only for large time horizons, therefore in the MPC open-loop problem, one can consider to fix the current velocity field $\mathbf{v}$ and solve only the linear convection-diffusion equation \eqref{eq:pde-heat-equation}. 

However, we will still carry out occasional simulations of the Boussinesq approximation in order to generate the source for so-called snapshots of the system state. From those we can then construct a local (with respect to time) reduced-order model of the convection-diffusion equation \eqref{eq:pde-heat-equation} system using POD. This reduced model summarizes information about the essential dynamics currently going on in the system and is much smaller than the original system which allows us to use it for fast computation of an approximated solution of optimal control problems within the MPC algorithm.

One may ask, why not performing POD directly for the coupled system \eqref{eq:Bous}, as in \cite{Rav2000,Rav2011}. Although this is a viable approach, the use of the velocity field $\textbf{v}$ as datum of the problem and not as unknown improves and simplifies further the POD approximation. This has the double advantage of using less basis functions and easier (mathematically and computationally) a-posteriori error estimator for the problem, which may result in tighter estimates compared to the one for the fully coupled system.

The last challenge is caused by the state constraints. These do not
guarantee the existence of regular ($L^2$) Lagrange multipliers
(cf. e.g.\ \cite{Tro13}) , introducing additional computational
effort. To overcome this issue, we employ a virtual control approach
together with a primal-dual active set strategy.

\section{Virtual control approach and PDASS for the linear-quadratic problem}
\label{sec:pdass}

As mentioned in Section~\ref{sec:setting}, as first simplification we assume that the velocity field $\textbf{v}$ is given and we focus on a model based only on the heat equation \eqref{eq:pde-heat-equation}-\eqref{eq:pde-heat-equation-initguess} with convection together with the boundary conditions \eqref{eq:pde-heat-equation-outside}-\eqref{eq:pde-heat-equation-control}. It is well-known that this equation admits a unique weak solution $y\in \Y=W(0,T;V)$ for any given control $u\in\U$; cf., e.g., \cite{DL00}. Let $\mathcal S: \U\to\Y$ be the control-to-state map for the heat equation, then we define the admissible set of controls 
\[
\Uad= \big\{u\in\U\,\big|\,u\text{ and }y=\mathcal S(u)\text{ satisfy \eqref{control_constraints} and \eqref{state_constraints}, respectively}\big\}.
\] 
We are, thus, interested in solving the strictly convex optimal control problem
\begin{equation}
\tag{$\mathbf{\hat P}$}
\label{eq:opt_cont_prob_in_u}
\min \mathcal{J}(u)\quad\text{s.t.}\quad u\in\Uad,
\end{equation}
which admits a unique solution $\bar u$ if the set $\Uad$ has non-empty interior; see, e.g., \cite{HPUU2009}.

Note that the difference between \eqref{eq:opt_cont_prob_in_u_bous}
and \eqref{eq:opt_cont_prob_in_u} relies on the different
control-to-state map, i.e. on the different underlying state model. As
said, we expect that $\bar u_\mathsf{B}$ and $\bar u$ will be close
for short time horizon $T$, since the buoyancy effect influences the temperature trajectory $y$ only on larger time scales. This will allow computing a good enough approximation of the velocity field $\textbf{v}$ for \eqref{eq:pde-heat-equation}, to be used as datum of the problem. We will explain in Section~\ref{sec:parallel} how we compute this approximation. 
What is important to clarify now is how to compute the solution to \eqref{eq:opt_cont_prob_in_u}.  As mentioned, the presence of pointwise state constraints leads to non-regular Lagrange multipliers; cf. \cite{Cas97,Ray97,RMRT08}. To overcome this issue, we apply a virtual control approach as in \cite{KR09,Mec19,MV18}. Let $w\in\W= L^2(Q)$ be an additional (artificial) control. For given $\varepsilon>0$ we relax the bilateral state constraints \eqref{state_constraints} as follows:
\begin{equation}
\label{relaxed_state_constraints}
\ya \leq \mathcal{S}(u) + \varepsilon w \leq \yb\quad\text{a.e. in } Q
\end{equation}
and introduce the artificial constraints
\begin{equation}
\label{artificial_constraints}
\wa \leq w(t,\bx) \leq \wb\quad\text{a.e. in } Q
\end{equation}
with fixed (safeguard) scalars $\wa\ll0$ and $\wb\gg0$, i.e., sufficiently small and large, respectively. Furthermore, we define the admissible set of controls
\[
\Zadeps = \big\{(u,w)\in\U\times\W\,\big|\,u, w\text{ and }\mathcal{S}(u)\text{ satisfy } \eqref{control_constraints}, \eqref{relaxed_state_constraints}, \eqref{artificial_constraints} \big\}.
\]
Now we are interested in solving
\begin{equation}
\tag{$\mathbf{\hat P}^\varepsilon$}
\label{eq:optimal_control_in_u_w}
\min J(u,w)\quad\text{s.t.}\quad(u,w)\in\Zadeps,
\end{equation}
where
\[
J(u,w)= \mathcal{J}(u)+\frac{\sigma}{2}\,{\|w\|}^2_\W
\]
with $\sigma>0$. Note that \eqref{eq:optimal_control_in_u_w} admits a unique solution $(\bar u^\varepsilon,\bar w^\varepsilon)\in \U\times\W$ for each $\varepsilon>0$ provided $\Zadeps$ has non-empty interior (cf. \cite{HPUU2009,Mec19}).

\begin{remark}
	\em
	We point out that the artificial bounds $\wa$, $\wb$ are necessary to guarantee $\bar u^\varepsilon \to\bar u$ in $\U$ as $\varepsilon\to 0$ with an order of convergence $O(\sqrt{\varepsilon})$; cf. \cite[Section~1.3.3]{Mec19}. They can be omitted if $\varepsilon$ is taken small but fixed. In this case, the existence of a unique solution is still guaranteed \cite{MV18}. Moreover, we can choose $\wa$ and $\wb$ as large in absolute value as one desires, because they do not have any physical meaning \cite{Mec19}.\hfill$\Diamond$
\end{remark}

Let us define the non-symmetric continuous weakly-coercive bilinear form $a(t;\cdot\,,\cdot):V\times V\to\mathbb{R}$ as
\[
a(t; \varphi,\psi) = \int_\Omega  {\textstyle\frac{\kappa}{\rho c_p}}\nabla\varphi\cdot\nabla\psi + (\textbf{v}(t,\cdot)\cdot \nabla\varphi)\psi\,\mathrm{d}\bx + \gamma\int_{\Gamma\setminus\Gamma_c} \varphi\psi\,\mathrm{d}\bs + \gamma_c\int_{\Gamma_c} \varphi\psi\,\mathrm{d}\bs.
\]
Then the first-order necessary optimality conditions for \eqref{eq:optimal_control_in_u_w} are as follows (cf. \cite{Mec19}):
\begin{theorem}
	\label{TheoremFOC}
	Let $\Zadeps$ have a non-empty interior. Suppose that the pair $(\bar u^\varepsilon,\bar w^\varepsilon)\in\Zadeps$ is the solution to \eqref{eq:optimal_control_in_u_w} with associated optimal state $\bar y^\varepsilon=\mathcal S(\bar u^\varepsilon)$, i.e. satisfying
	\begin{subequations}
		\label{VC:OptSyst}
		\begin{equation}
		\label{StateEquation}
		\begin{aligned}
		\frac{\mathrm d}{\mathrm dt}\,{\langle \bar y^\varepsilon(t),\varphi\rangle}_H+a(t;\bar y^\varepsilon(t),\varphi)&=\gamma y_\text{out}(t)\int_{\Gamma\setminus\Gamma_c} \varphi\,\mathrm{d}\bs + \gamma_c u(t)\int_{\Gamma_c} \varphi\,\mathrm d \bs,\\
		\bar y^\varepsilon(0)&= y_\circ
		\end{aligned}
		\end{equation}
		for all $\varphi\in V$ and a.e. in $[0,T]$, where $H= L^2(\Omega)$. Then, there exist unique Lagrange multipliers $\bar p^\varepsilon\in \Y$, $\bar\beta^\varepsilon,\bar \vartheta^\varepsilon \in\W$ and $\bar\alpha^\varepsilon\in\U$ satisfying the dual equation
		\begin{equation}
		\label{DualEquation}
		\begin{aligned}
		-\frac{\mathrm d}{\mathrm dt}\,{\langle \bar p^\varepsilon(t),\varphi\rangle}_H+a(t;\varphi,\bar p^\varepsilon(t))&=-{\langle\bar\beta^\varepsilon(t),\varphi\rangle}_H,\\
		\bar p^\varepsilon(T)&=-\bar \beta^\varepsilon(T)
		\end{aligned}
		\end{equation}
		for all $\varphi\in V$ and a.e. in $[0,T]$ and the optimality system
		\begin{align}
		\label{OptConda} \bar u^\varepsilon -\gamma_c \int_{\Gamma_c}\bar p^\varepsilon\mathrm{d} s +\bar\alpha^\varepsilon &=0&&\text{in }\U,\\
		\label{OptCondb} \sigma \bar w^\varepsilon+\varepsilon\bar\beta^\varepsilon+\bar\vartheta^\varepsilon&=0&&\text{in }\W.
		\end{align}
		Moreover,
		\begin{align}
		\label{NCP-1}
		\bar\beta^\varepsilon&=\max\big\{0,\bar\beta^\varepsilon+\eta(\bar y^\varepsilon+\varepsilon\bar w^\varepsilon-\yb)\big\}+\min\big\{0,\bar\beta^\varepsilon+\eta(\bar y^\varepsilon+\varepsilon\bar w^\varepsilon-\ya\big\},\\
		\label{NCP-2}
		\bar\alpha^\varepsilon&=\max\big\{0,\bar\alpha^\varepsilon+\eta_u(\bar u^\varepsilon-\ub)\big\}+\min\big\{0,\bar\alpha^\varepsilon+\eta_u(\bar u^\varepsilon-\ub)\big\}, \\
		\label{NCP-3}
		\bar\vartheta^\varepsilon&=\max\big\{0,\bar\vartheta^\varepsilon+\eta_w(\bar w^\varepsilon-\wb)\big\}+\min\big\{0,\bar\vartheta^\varepsilon+\eta_w(\bar w^\varepsilon-\wa)\big\}
		\end{align}
	\end{subequations}
	for arbitrarily chosen $\eta,\eta_u,\eta_w>0$, where the max- and min-operations are interpreted componentwise in the pointwise everywhere sense.
\end{theorem}

\begin{remark}
	\em
	The well-defined bounded solution operator $\mathcal A:\W\to\Y$ is defined as follows: for given $\beta\in\W$ the function $p=\mathcal A\beta$ is the unique solution to
	\begin{align*}
	-\frac{\mathrm d}{\mathrm dt}\,{\langle p(t),\varphi\rangle}_H+a(t;\varphi,p(t))&=-{\langle \beta(t),\varphi\rangle}_H&&\forall\varphi\in V\text{ a.e. in }[0,T),\\
	p(T)&=-\beta(T)&&
	\end{align*}
	for given $\beta\in\W$; cf. Remark~2.3 in \cite{MV18}. Then, $\bar p^\varepsilon=\mathcal A(\bar\beta^\varepsilon)$ solves \eqref{DualEquation}.\hfill$\Diamond$
\end{remark}

Letting $\nu=(u,w,\theta)\in \mathcal N=\U\times\W\times\W$, we define the active sets corresponding to \eqref{VC:OptSyst} as
\begin{subequations}
	\label{ActInactSets}
	\begin{equation}
	\label{ActiveSets}
	\begin{aligned}
	\AaU(\nu)&=\big\{t\in(0,T)\,\big|\,\alpha(\nu)+u-\ua<0\text{ a.e.}\big\},\\
	\AbU(\nu)&=\big\{t\in(0,T)\,\big|\,\alpha(\nu)+u-\ub>0\text{ a.e.}\big\},\\
	\AaW(\nu)&=\bigg\{(t,\bx)\in Q\,\big|\,\beta(\nu)+\frac{\sigma}{\varepsilon^2}\big(y(\nu)+\varepsilon w-\ya\big)<0\text{ a.e.}\bigg\},\\
	\AbW(\nu)&=\bigg\{(t,\bx)\in Q\,\big|\,\beta(\nu)+\frac{\sigma}{\varepsilon^2}\big(y(\nu)+\varepsilon w-\yb\big)>0\text{ a.e.}\bigg\}, \\
	\AaVT (\nu)& = \big\{(t,\bx)\in Q\,\big|\,w-\wa<0\text{ a.e.}\big\}, \\
	\AbVT (\nu)& = \big\{(t,\bx)\in Q\,\big|\,w-\wb>0\text{ a.e.}\big\}. 
	\end{aligned}
	\end{equation}
	Similarly, the associated inactive sets are
	\begin{equation}
	\label{InactiveSets}
	\begin{aligned}
	\IU(\nu)&=(0,T)\setminus\big(\AaU(\nu)\cup\AbU(\nu)\big),\\
	\IW(\nu)&=Q\setminus\big(\AaW(\nu)\cup\AbW(\nu)\big),\quad\IVT(\nu)= Q\setminus\big(\AaVT(\nu)\cup\AbVT(\nu)\big).
	\end{aligned}
	\end{equation}
\end{subequations}
One can solve \eqref{VC:OptSyst} applying a primal-dual active set strategy or, equivalently, a semi-smooth Newton method; see, e.g., \cite{HIK02,IK03,Mec19}. Doing so, at each iteration $k$ of the semi-smooth Newton method we have to solve the following system
\begin{subequations}
	\label{PDASS-System}
	\begin{align}
	\label{PDASS-System-a}
	\gamma_c\int_{\Gamma_c}p^{k+1}\,\mathrm{d}\bs-u^{k+1}&=0 &&\text{in }\IU(\nu^k),\\
	\label{PDASS-System-b}
	u^{k+1}&=\ua&&\text{in }\AaU(\nu^k),\\
	\label{PDASS-System-c}
	u^{k+1}&=\ub&&\text{in }\AbU(\nu^k),\\
	\label{PDASS-System-d}
	w^{k+1}&=0&&\text{in }\IW(\nu^k),\\
	\label{PDASS-System-e}
	y^{k+1}+\varepsilon\,w^{k+1}&=\ya&&\text{in }\AaW(\nu^k),\\
	\label{PDASS-System-f}
	y^{k+1}+\varepsilon\,w^{k+1}&=\yb&&\text{in }\AbW(\nu^k), \\
	\label{PDASS-System-g}
	\theta^{k+1}&=0 && \text{in }\IVT(\nu^k), \\
	\label{PDASS-System-h}
	\sigma w^{k+1}-\varepsilon \theta^{k+1} &= \sigma \wa && \text{in }\AaVT(\nu^k), \\
	\label{PDASS-System-i}
	\sigma w^{k+1}-\varepsilon \theta^{k+1} &= \sigma \wb && \text{in }\AbVT(\nu^k).  
	\end{align}
\end{subequations}
Taking into account \eqref{VC:OptSyst}, conditions \eqref{PDASS-System} lead to the following coupled equations:
\begin{subequations}
	\label{eq:pdass_coupled_system}
	\begin{align}
	&\frac{\mathrm d}{\mathrm dt} {\langle y^{k+1}(t),\varphi \rangle}_H+a(y^{k+1}(t),\varphi) -\gamma_c \mathcal{H}^{k}(t; p^{k+1}(t))\int_{\Gamma_c}\varphi\,\mathrm{d}\bs\\\nonumber
	&\hspace{0.36\textwidth}=\gamma y_\text{out}(t)\int_{\Gamma\setminus\Gamma_c} \varphi\,\mathrm{d}\bs+  \gamma_cr^k(t)\int_{\Gamma_c} \varphi\,\mathrm{d}\bs, \\
	&\label{eq:pdass_dual} -\frac{\mathrm d}{\mathrm dt}\,{\langle p^{k+1}(t),\varphi\rangle}_H+a(t;\varphi,p^{k+1}(t))  +\frac{\sigma}{\varepsilon^2}\,\left\langle \mathcal{G}^{k,\varepsilon}(t;y^{k+1}(t)),\varphi\right\rangle_H
	\\\nonumber & \hspace{0.63\textwidth} = \frac{\sigma}{\varepsilon^2}\,{\langle r^{k,\varepsilon}(t),\varphi\rangle}_H 
	\end{align} 
	for all $\varphi\in V$ and a.e. in $(0,T)$, together with the conditions
	\begin{align}
	y^{k+1}(0) &=y_\circ, \\
	\label{pT}
	p^{k+1}(T) + \frac{\sigma}{\varepsilon^2}\,\mathcal{G}^{k,\varepsilon}(T;y^{k+1}(T)) & = \frac{\sigma}{\varepsilon^2}\, r^{k,\varepsilon}(T),
	\end{align}
\end{subequations}
where
\begin{align*}
\mathcal{G}^{k,\varepsilon}(t;y^{k+1}(t))&= y^{k+1}(t)\sum_{i=1}^6 \chi_{\mathscr{A}_i(\nu^k)}(t)+ \frac{1}{\varepsilon}y^{k+1}(t)\sum_{i=3}^6 \chi_{\mathscr{A}_i(\nu^k)}(t) \\
\mathcal{H}^k(t;p^{k+1}(t)) &= \chi_{\IU(z^k)}(t)\gamma_c\int_{\Gamma_c} p^{k+1}(t)\mathrm d s \\
r^{k}(t) &= \chi_{\AaU(z^k)}(t)\ua(t)+\chi_{\AbU(z^k)}(t)\ub(t) \\
r^{k,\varepsilon}(t) &= \ya\chi_{\AWon(\nu^k)}(t) + \yb\chi_{\AWtw(\nu^k)}(t) - \wa\big(\chi_{\AWth(\nu^k)}(t)+\chi_{\AWfi(\nu^k)}(t)\big)\\ & \quad - \wb\big(\chi_{\AWfo(\nu^k)}(t)+\chi_{\AWsi(\nu^k)}(t)\big) \\ & \quad + \frac{\varepsilon+1}{\varepsilon} \ya \big( \chi_{\AWth(\nu^k)}(t)+\chi_{\AWfo(\nu^k)}(t)\big) \\
& \quad + \frac{\varepsilon+1}{\varepsilon} \yb \big( \chi_{\AWfi(\nu^k)}(t)+\chi_{\AWsi(\nu^k)}(t)\big),
\end{align*}
with $\chi_\bullet(t)$ being the indicator functions of the sets
\begin{align*}
\AWon(\nu^k)&= \AaW(\nu^k)\cap\IVT(\nu^k), & 
\AWtw(\nu^k)&=\AbW(\nu^k)\cap\IVT(\nu^k),\\
\AWth(\nu^k)&=\AaW(\nu^k)\cap\AaVT(\nu^k), &
\AWfo(\nu^k)&=\AaW(\nu^k)\cap \AbVT(\nu^k), \\
\AWfi(\nu^k)&=\AbW(\nu^k)\cap \AaVT(\nu^k),  &
\AWsi(\nu^k)&=\AbW(\nu^k)\cap \AbVT(\nu^k),
\end{align*}
respectively. Summarizing, the coupled system of equations \eqref{eq:pdass_coupled_system} can be formulated only in the variables $y^{k+1}$ and $p^{k+1}$:
\begin{equation}
\label{OpSystem}
\left(
\begin{array}{cc}
\mathcal A_{11}^k&\mathcal A_{12}^k\\[1ex]
\mathcal A_{21}^k&\mathcal A_{22}^k
\end{array}
\right)\left(
\begin{array}{c}
y^{k+1}\\[1ex]
p^{k+1}
\end{array}
\right)=\left(
\begin{array}{c}
\mathcal Q_1(\nu^k;\ua,\ub,\gamma_c,\gamma,y_\text{out})\\[1ex]
\mathcal Q_2(\nu^k;\ya,\yb,\wa,\wb,\varepsilon,\sigma)
\end{array}
\right),
\end{equation}
where the operators $\mathcal A_{11}^k$ and $\mathcal A_{22}^k$ are of the form
\[
\mathcal A_{11}^k=\mathcal C+\tilde{\mathcal A}_{11}^k, \quad\mathcal A_{22}^k=\mathcal C^\star+\tilde{\mathcal A}_{22}^k
\]
and $\mathcal C$ stands for the heat-convection operator, which does not depend on $k$. A discretization of \eqref{OpSystem} leads to a discretized system of the form
\begin{equation}
\label{DiscOpSystem}
\left(
\begin{array}{cc}
\mathrm A_{11}^k&\mathrm A_{12}^k\\[1ex]
\mathrm A_{21}^k&\mathrm A_{22}^k
\end{array}
\right)\left(
\begin{array}{c}
\mathrm y^{k+1}\\[1ex]
\mathrm p^{k+1}
\end{array}
\right)=\left(
\begin{array}{c}
\mathrm Q_1^k\\[1ex]
\mathrm Q_2^k
\end{array}
\right)
\end{equation}
with
\[
\mathrm A_{11}^k=\mathrm C+\tilde{\mathrm A}_{11}^k, \quad\mathrm A_{22}^k=\mathrm C^\top+\tilde{\mathrm A}_{22}^k,
\]
where $\mathrm C$ stands for the discretized heat-convection operator, which is again independent of $k$. For further details, we refer to \cite{Mec19}. We summarize the PDASS in Algorithm~\ref{Alg:PDASS}.
\begin{algorithm}
	\caption{(Primal-dual active set strategy)}
	\label{Alg:PDASS}
	\begin{algorithmic}[1]
		\STATE Choose starting value $\nu^0=(u^0,w^0,\vartheta^0)\in\N$, set $k=0$ and {\texttt flag}~= {\texttt false};
		\STATE Determine $y^0=\mathcal Su^0$ and $p^0=- \mathcal A(\sigma w^0+\vartheta^0)/\varepsilon$.
		\STATE Determine $\AaU(\nu^0)$, $\AbU(\nu^0)$, $\IU(\nu^0)$ for $i=1,\ldots,m$, $\AaW(\nu^0)$, $\AbW(\nu^0)$, $\IW(\nu^0)$ and $\AaVT(\nu^0)$, $\AaVT(\nu^0)$ and $\IVT(\nu^0)$ from \eqref{ActInactSets};
		\REPEAT
		\STATE Compute the solution $(y^{k+1},p^{k+1})$ by solving \eqref{DiscOpSystem};
		\STATE Compute $\nu^{k+1}=(u^{k+1},w^{k+1},\vartheta^{k+1})\in\N$ from \eqref{PDASS-System} and set $k=k+1$;
		\STATE Determine $\AaU(\nu^k)$, $\AbU(\nu^k)$, $\IU(\nu^k)$ for $i=1,\ldots,m$, $\AaW(\nu^k)$, $\AbW(\nu^k)$, $\IW(\nu^k)$ and $\AaVT(\nu^k)$, $\AaVT(\nu^k)$ and $\IVT(\nu^k)$ from \eqref{ActInactSets};
		\IF{$\mathscr A_{{\mathsf a}i}^\U(\nu^k)=\mathscr A_{{\mathsf a}i}^\U(\nu^{k-1})$ {\textbf and} $\mathscr A_{{\mathsf b}i}^\U(\nu^k)=\mathscr A_{{\mathsf b}i}^\U(\nu^{k-1})$ for $i=1,\dots,m$}
		\IF{$\AaW(\nu^k)=\AaW(\nu^{k-1})$ {\textbf and} $\AbW(\nu^k)=\AbW(\nu^{k-1})$}
		\IF{$\AaVT(\nu^k)=\AaVT(\nu^{k-1})$ {\textbf and} $\AbVT(\nu^k)=\AbVT(\nu^{k-1})$}
		\STATE Set {\texttt flag}~=~{\texttt true};
		\ENDIF
		\ENDIF
		\ENDIF
		\UNTIL{{\texttt flag}~=~{\texttt true};}
	\end{algorithmic}
\end{algorithm}
\begin{remark}
	\label{Remark:PDASS}
	\em
	The discrete linear system \eqref{DiscOpSystem} can be obtained discretizing the primal and dual equations \eqref{eq:pdass_coupled_system}, for example, with piecewise linear finite elements (FE) in space and the implicit Euler scheme in time. Let $V_h\subset V$ be the FE space with $\dim V_h = N_\bx$ and $N_t$ be the number of time steps. Notice that the matrix in \eqref{DiscOpSystem} has the size $(2N_tN_\bx)\times(2N_tN_\bx)$. Clearly, its dimension may cause problems regarding memory consumption and computational time. In our application, in fact, $N_\bx$ is not small, since a coarse FE grid can not capture the complexity of the dynamics involved. Therefore, we apply POD-based reduced-order modeling to approximate the solution of \eqref{DiscOpSystem} by a reduced-order system of dimension $2\ell N_t$ with $\ell\ll N_\bx$, gaining computational time but paying in approximation; cf. Section~\ref{sec:POD}. The number of time steps may not be small either, since the time horizon $T$ could be extremely large (e.g. one month) and the time step small (e.g. a minute). This is not a disadvantage for us, since because we apply an MPC method, we iteratively solve optimal control subproblems on a smaller time horizon, computing a feedback control that approximates the optimal control $\bar u^\varepsilon$; cf.\ Section~{\ref{sec:MPC}}. Finally, note that the matrix in \eqref{DiscOpSystem} has a sparse block structure that allows to store only the FE matrices, e.g., the mass and stiffness matrices, which have dimension $N_\bx\times N_\bx$. These blocks will then be projected into the POD space. Doing that we can save memory significantly, in particular when the number of FE nodes is  large. Note that, in comparison to the Boussinesq approximation, we gain even more, since in general computing the solution of this model requires an additional picewise quadratic FE approximation of the velocity space, which would definitely become infeasible in combination with the PDASS.
\end{remark}

\section{POD and a-posteriori error estimate}
\label{sec:POD}
Model order reduction and, in particular, POD is a well established field of numerical analysis, which was intensively developed in the past 20 years. The key idea is to construct a low-dimensional subspace $V_\ell\subset V_h=\mathrm{span}\,\{\varphi_1,\ldots,\varphi_{N_\bx}\}$, spanned by the so-called POD basis functions. For the sake of completeness, in the first part of this section, we describe briefly the POD method in the discrete case, but we focus on one snapshot ensemble for brevity. For more details we refer the reader to \cite{GV17}, for instance. Let $t_j=j\Delta t$, $j=0,\ldots,n$, be the time discretization for a suitable time interval with constant time step $\Delta t>0$. By $\{y_j\}_{j=0}^n\subset V_h$ we denote a given discrete trajectory, also called snapshots ensemble. We define the snapshot subspace $\mathscr V_n= \mathrm{span}\,\{y_0,\ldots,y_n\}\subset V_h$. Let $d_n = \dim\mathscr V_n\le n+1$. For any $\ell\in\{1,\ldots,d_n\}$ we are interested in finding an orthonormal set $\{\psi_i\}_{i=1}^\ell\subset\mathscr V_n$, which is the solution of the following minimization problem
\begin{equation}
\label{eq:POD_prob}
\min \sum_{j=0}^n \alpha_j \Big\|y_j-\sum_{i=1}^\ell {\langle y_j,\psi_i\rangle}_V\,\psi_i\Big\|_V^2\quad\text{s.t.}\quad{\langle\psi_i,\psi_j\rangle}_V=\delta_{ij},
\end{equation}
where $\delta_{ij}$ is the Kronecker delta and $\alpha_j>0$ are given (trapezoidal) weights. We set $V^\ell=\mathrm{span}\,\{\psi_1,\ldots,\psi_\ell\}\subset V_h$. Obviously, if $\dim V_\ell=\ell\ll m$ holds we gain a computational speed-up. From $y_j\in V_h$ ($0\le j\le n$) we infer that there exists a snapshot matrix $Y\in\mathbb R^{N_\bx\times(n+1)}$, such that
\begin{align*}
y_j(\bx)=\sum_{i=1}^{N_\bx}Y_{ij}\varphi_i(\bx)\quad\text{for }\bx\in\Omega.
\end{align*}
Introducing the weighting matrices $D=\mathrm{diag}\,(\alpha_0,\ldots,\alpha_n)\in\mathbb{R}^{(n+1)\times (n+1)}$ and $W=((W_{ij}))\in\mathbb{R}^{N_\bx\times N_\bx}$ with $W_{ij}=\langle\varphi_j,\varphi_i\rangle_V$, a solution $\{\psi_j\}_{j=1}^\ell$ to \eqref{eq:POD_prob} can be computed as follows:
\begin{itemize}
	\item Solve the symmetric eigenvalue problem
	\begin{align*}
	\big(D^{1/2}Y^\top WYD^{1/2}\big)\bm\phi_j = \lambda_j\bm\phi_j\quad\text{for }j=1,\ldots,\ell
	\end{align*}
	with $\lambda_1\ge\ldots\ge \lambda_\ell>0$.
	\item Define the vectors $\bm\psi_j=YD^{1/2}\bm\phi_j/\sqrt{\lambda_j}\in\mathbb R^{N_\bx}$ for $1\leq j\leq \ell$ and the matrix $\Psi=[\bm\psi_1|\ldots|\bm\psi_\ell]\in\mathbb R^{N_\bx\times\ell}$.
	\item Set
	\begin{align*}
	\psi_j(\bx)=\sum_{i=1}^{N_\bx}\Psi_{ij}\varphi_i(\bx)\quad\text{for }j=1,\ldots,\ell\text{ and }\bx\in\Omega.
	\end{align*}
\end{itemize}

Using a Galerkin ansatz, the reduced-order version of \eqref{eq:pdass_coupled_system} can be derived replacing the test functions $\varphi$ by the POD basis $\left\{\psi_i\right\}_{i=1}^\ell$. We indicate with $\mathcal{S}^\ell$ the POD solution operator corresponding to the reduced-order model. We can then construct the discretized system
\begin{equation}
\label{DiscOpSystemPOD}
\left(
\begin{array}{cc}
\mathrm A_{11}^{k,\ell}&\mathrm A_{12}^{k,\ell}\\[1ex]
\mathrm A_{21}^{k,\ell}&\mathrm A_{22}^{k,\ell}
\end{array}
\right)\left(
\begin{array}{c}
\mathrm y^{k+1,\ell}\\[1ex]
\mathrm p^{k+1,\ell}
\end{array}
\right)=\left(
\begin{array}{c}
\mathrm Q_1^{k,\ell}\\[1ex]
\mathrm Q_2^{k,\ell}
\end{array}
\right),
\end{equation}
of which the matrix belongs to $\mathbb{R}^{(2\ell N_t)\times (2\ell N_t)}$. From \eqref{DiscOpSystemPOD}, one can easily implement a POD-based PDASS following the structure of Algorithm~\ref{Alg:PDASS}. Now, a first question is how to check that the chosen number of POD basis $\ell$ is good enough to guarantee a small approximation error in reconstructing the snapshots. For this, the following error formula holds
\begin{equation}
\label{aprioriellselection}
\sum_{j=0}^n \alpha_j \bigg\|y_j-\sum_{i=1}^\ell {\langle y_j,\psi_i\rangle}_V \psi_i\bigg\|_V^2 = \sum_{i=\ell+1}^{d_n}\lambda_i.
\end{equation}
More importantly, to check that the constructed POD basis is able to reconstruct the full-order optimal solution $\bar u^\varepsilon$, we need an a-posteriori error estimator. Let $\bar u^{\varepsilon,\ell}$ be the solution of the POD-based PDASS, then the following result from \cite{Mec19} holds:

\begin{proposition}
	\label{prop:err_est}
	Let $(\bar u^\varepsilon,\bar w^\varepsilon)\in\Zadeps$ be the optimal solution to \eqref{eq:optimal_control_in_u_w}. Suppose that $z^\mathsf{ap}=(u^\mathsf{ap},w^\mathsf{ap})\in\Zadeps$ is given arbitrarily. Then, there exists a perturbation $\zeta=(\zeta^u,\zeta^w)\in\U\times\W$, which is independent of $(\bar u^\varepsilon,\bar w^\varepsilon)$, so that
	\begin{equation}
	\label{APostError}
	\|\bar u^\varepsilon-u^\mathsf{ap}\|_\U + \|\bar w^\varepsilon-w^\mathsf{ap}\|_\W\le\frac{1}{\sigma^\mathsf{ap}}\,{\|\mathcal T^\star\zeta\|}_\U,
	\end{equation}
	where $\sigma^\mathsf{ap}:=\min\{\sigma,1\}>0$ and $\mathcal T^\star$ is defined as
	\[
	\mathcal T^\star=\left(
	\begin{array}{cc}
	\mathcal I_\U& \mathcal S^\star \\
	0 & \varepsilon \mathcal I_\W 
	\end{array}
	\right):\U\times\W\to\U\times\W,
	\]
	where $\mathcal I_\U$ and $\mathcal I_\W$ are the identities on $\U$ and $\W$, respectively. The perturbation $\zeta$ is computed as follows: Let $\xi=(\xi^u,\xi^w)\in\U\times\W$ be given as the solution of the linear system $\mathcal T^\star\xi=\nabla\hat J(z^\mathsf{ap})$, i.e.,
	\begin{equation}
	\label{SystemAPostError}
	\left(\begin{array}{cc}
	\mathcal I_\U& \mathcal S^\star \\
	0&\varepsilon \mathcal I_\W
	\end{array}\right)
	\left(\begin{array}{c}
	\xi^u\\
	\xi^w
	\end{array}\right)
	=\left(
	\begin{array}{c}
	u^\mathsf{ap}\\
	\sigma w^\mathsf{ap}
	\end{array}
	\right),
	\end{equation}
	then
	\begin{subequations}
		\label{PertZeta}
		\begin{equation}
		\label{PertZeta-a}
		\zeta^u(t)=\left\{\begin{array}{ll}
		-\min\{0,\xi^u(t)\}&\text{for }t\in \AaU(z^\mathsf{ap}),\\[1mm]
		-\max\{0,\xi^u(t)\}&\text{for }t\in\AbU(z^\mathsf{ap}),\\[1mm]
		-\xi^u(t)&\text{for }t\in\IU(z^\mathsf{ap})
		\end{array}\right.
		\end{equation}
		and
		\begin{equation}
		\label{PertZeta-b}
		\zeta^w(t,\bx)=\left\{
		\begin{array}{ll}
		-\min\{0,\xi^w(t,\bx)\}&\text{for }(t,\bx)\in\AaW(z^\mathsf{ap}),\\[1mm]
		-\max\{0,\xi^w(t,\bx)\}&\text{for }(t,\bx)\in\AbW(z^\mathsf{ap}),\\[1mm]
		-\xi^w(t,\bx)&\text{for }(t,\bx)\in\IW(z^\mathsf{ap}).
		\end{array}\right.
		\end{equation}
	\end{subequations}
\end{proposition}
\begin{remark}
	\em
	Note that we can easily decouple the equations in system \eqref{SystemAPostError} by computing $\xi^w=\sigma_w w^\mathsf{ap}/\varepsilon$ from the second equation. Then $\xi^u = -\mathcal S^\star \xi^w + u^\mathsf{ap} = -\sigma\mathcal S^\star w^\mathsf{ap}/\varepsilon + u^\mathsf{ap}$. Moreover, we recall that $\mathcal S^\star: W_0(0,T)'\to\U$ denotes the dual operator of the linear solution operator $\mathcal S$; see \cite[Lemma 2.4]{TV09}. Finally, we remark that in our numerical realization $z^\mathsf{ap}$ is given by the POD suboptimal solution pair $(\bar u^{\varepsilon,\ell},\bar w^{\varepsilon,\ell})$. Thus, \eqref{APostError} can be utilized as an a-posteriori error estimate. For further details regarding the a-posteriori error estimator, we refer to \cite{Mec19}.\hfill$\Diamond$
\end{remark} 

\section{Economic MPC}
\label{sec:MPC}

MPC is a well-established method for computing a closed-loop control for a dynamical system on an infinite (or large) time horizon. 
MPC splits the solution of such problems into the consecutive solution of problems on a small finite time horizon, the so called prediction horizon.
The goal is to reduce the complexity of the problem in time, while at the same time being able to react to model uncertainties, disturbances, and possible parameter changes. We are interested in economic MPC \cite{AmRA11,FaGM18}. This term comprises all MPC schemes in which the cost function does not merely penalize the distance to a pre-defined target. Here we do not have such a pre-defined target. Rather, the goal is to minimize the control effort, while satisfying certain control and state constraints; see, e.g., \cite{GP17,GP20,Pir20}. Applying MPC to the virtual control problem \eqref{eq:optimal_control_in_u_w} has several advantages, which will be more clear after having introduced the method. Let $N\Delta t$, $N\in\N$, be the prediction horizon. We define the Hilbert spaces
\[
\UNj = L^2(t_j,t_{j+N}), \quad \WNj = L^2((t_j,t_{j+N})\times\Omega).
\]
for $j=0,\ldots,N_t$ and with the notation $t_j= j\Delta t$ for $j>N_t$. Then the MPC cost functional is
\[
J^j_N(u,w)= \frac{1}{2} \|u\|^2_{\UNj} + \frac{\sigma}{2}\|w\|^2_{\WNj}
\] 
and the MPC admissible sets are
\[
\begin{aligned}
\ZadepsNj = \big\{& (u,w)\in\UNj\times\WNj: u, w, \mathcal{S}u+\varepsilon w \text{ satisfy }  \\ &\ua(t)\leq u(t)\leq \ub(t) \text{ a.e. in } (t_j,t_{j+N}], \\ & \wa\leq w(t,\bx)\leq \wb, \, \ya\leq \mathcal{S}u+\varepsilon w\leq \yb \text{ a.e. in } (t_j,t_{j+N}]\times\Omega, \\ & \text{ respectively} \big\},
\end{aligned}
\]
for $j=0,\ldots,N_t$. Furthermore, note that the fully discretized version of \eqref{StateEquation} can be rewritten as the following discrete dynamical system
\begin{equation}
\label{eq:discrete_dynamical_system}
\begin{aligned}
y_{j+1} = f(j,y_j,u_j)\text{ for } j=0,\ldots,N_t,\quad y_0 = y_\circ,
\end{aligned}
\end{equation}
where the function $f$ comprises all the formulas that come from the discretization and the model, which we avoid to write for brevity. The MPC method is then summarized in Algorithm~\ref{Alg:MPC}.
\begin{algorithm}
	\caption{Model predictive control \label{Alg:MPC}}
	\begin{algorithmic}[1]
		\STATE Choose an initial state $y_0$, an MPC prediction horizon $N$ and the regularization parameter $\varepsilon>0$;
		\FOR {each time instant $j=0,1,2,\ldots,N_t$} 
		\STATE \label{MPC_measurestep} Measure the current state of the system $y^j$ and the model parameters;
		\STATE Solve the optimal control problem 
		\begin{equation}
		\label{eq:opt_control_problem_MPC}
		\min J_{N}^j(u,w)\quad\text{s.t.}\quad(u,w)\in\ZadepsNj
		\end{equation}
		to obtain the open-loop optimal control sequence $\bar u^\varepsilon_N$;
		\STATE Apply the first element of $\bar u^\varepsilon_N$ as a control to the system \eqref{eq:discrete_dynamical_system} during the next sampling period, i.e.\ use the feedback law $\mu_N(y_j) = \bar u^\varepsilon_N(0)$.
		\ENDFOR
	\end{algorithmic}
\end{algorithm} 
To guarantee well-posedness and existence of the optimal control to each open-loop optimal control problem in Algorithm~\ref{Alg:MPC}, we assume that the admissible sets $\ZadepsNj$ have non-empty interior for all $j=0,\ldots,N_t$.

Algorithm~\ref{Alg:MPC} produces a sequence of optimal controls in
discrete time and uses the first element of each optimal control
sequence as a feedback control value. This results in a closed-loop
sub-optimal trajectory, which solves
\eqref{eq:discrete_dynamical_system} for $u_j= \mu_N(y_j)$, 
$j=0,\ldots,N_t$. Note that in step~\ref{MPC_measurestep} of the MPC
algorithm it is possible to measure (and update) the parameters of the
model. This is particularly advantageous for our application, since usually we do not have, e.g., a precise weather forecast for the entire long-time horizon $T$. In this way, we can update the outside temperature while the algorithm runs. To determine $y^j$ typically an estimator is
needed to process available sensor data properly to reconstruct the
spatial-temporal distribution. In addition, we can take advantage for the velocity field $\textbf{v}$. We mentioned already that this will be provided from an approximation of the Boussinesq velocity fields. We thus gain the possibility of updating this approximation while the simulation runs, improving the accuracy of our simplified model. 

The optimal control problems \eqref{eq:opt_control_problem_MPC} are solved on a much smaller time horizon $(t_j,t_{j+N})$ since $N\ll N_t$. This guarantees a good computational speed-up. For general optimal control problems it is not necessarily the case that the MPC closed-loop trajectory is approximately optimal on the long time horizon $N_t$. However, in recent years structural conditions were discovered that allow to conclude this property for $N$ large enough; cf.\ \cite[Chapter 7]{GrPa17} and \cite{FaGM18,GP20}. The most important among these conditions are the turnpike property and the closely related strict dissipativity property \cite{TreZ15,GruM16}. The turnpike property can be seen as a similarity condition for optimal trajectories on different time horizons and we will check numerically that this property holds for the optimal control problem considered in this paper; cf. Section~\ref{sec:turnpike}. It is also possible to check the quality of the MPC horizon and trajectory with a-posteriori error estimators, such as in \cite{GMPV20,GP09}, but we do not utilize them in this work.

What remains to be explained is how to successfully combine MPC and POD. The challenge is to guarantee a small approximation error for the reduced-order model while the MPC algorithm proceeds. This is clearly a big challenge if one thinks about the fact that the problem parameters are not only time-variant, but that the parameter values at a certain time instant in the future may also change when new measurement information becomes available. When this happens, it means that we need to update the POD model. Luckily, the a-posteriori error estimate from Proposition~\ref{prop:err_est} is a useful tool to check whether this should be done or not.
\begin{algorithm}[t]
	\caption{MPC-POD \label{Alg:MPCPOD}}
	\begin{algorithmic}[1]
		\STATE Choose an initial state $y_0$, an MPC prediction horizon $N$, the regularization parameter $\varepsilon>0$, and tolerances $\tau_1$ for the a-priori estimate \eqref{aprioriellselection} and $\tau_2$ for the a-posteriori estimate;
		\STATE Set $t_0=0$, $y_0(0)=y_\circ$, and {\texttt flag} = {\texttt true};
		\FOR{ $j=0,1,2,\dots,N_t$ }
		\STATE Measure the current state $y_j$ of the system at time $t_j$;
		\IF { {\texttt flag} = {\texttt true} }
		\STATE Update the model parameters;
		\STATE \label{Alg:MPCPOD:step:compute_snapshots} Compute primal and dual snapshots;
		\STATE \label{Alg:MPCPOD:step:compute_basis} Compute a POD basis $\{\psi_i\}_{i=1}^\ell$ of rank $\ell$ according to $\tau_1$;
		\STATE Set {\texttt flag} = {\texttt false};
		\ENDIF
		\STATE Solve the MPC-POD optimal control problem
		\begin{equation}
		\label{eq:opt_control_problem_MPCPOD}
		\min J^j_N(u,w) \text{ s.t. } (u,w)\in \ZadepsellNj
		\end{equation}
		to obtain the open-loop suboptimal control sequence $\bar u^{\varepsilon,\ell}_N$;
		\STATE \label{Alg:MPCPOD:step:compute_feedback} Apply the first element of $\bar u^{\varepsilon,\ell}_N$ as a control to the system \eqref{eq:discrete_dynamical_system} during the next sampling period, i.e.\ use the feedback law $\mu_N(y_j) = \bar u^{\varepsilon,\ell}_N(0)$;
		\STATE Compute the a-posteriori error estimate $e= {\|\mathcal T^\star\zeta\|}_{\U\times\W}$ from \eqref{APostError};
		\IF { $e> \tau_2$}
		\STATE Set {\texttt flag} = {\texttt true};
		\ENDIF
		\ENDFOR
	\end{algorithmic}
\end{algorithm}
We report the MPC-POD algorithm in Algorithm~\ref{Alg:MPCPOD}. For further details, we refer to \cite{Mec19,MV19}. The admissible sets in \eqref{eq:opt_control_problem_MPCPOD} are defined as
\[
\begin{aligned}
\ZadepsellNj = \big\{& (u,w)\in\UNj\times\WNj: u, w, \mathcal{S}^\ell u+\varepsilon w \text{ satisfy }  \\ &\ua(t)\leq u(t)\leq \ub(t) \text{ a.e. in } (t_j,t_{j+N}], \\ & \wa\leq w(t,\bx)\leq \wb, \, \ya\leq \mathcal{S}^\ell u+\varepsilon w\leq \yb \text{ a.e. in } (t_j,t_{j+N}]\times\Omega, \\ & \text{ respectively} \big\},
\end{aligned}
\]
where we recall that $\mathcal{S}^\ell$ is the POD solution operator for the reduced-order state equation. Obviously, to guarantee existence and uniqueness of the suboptimal controls solution to \eqref{eq:opt_control_problem_MPCPOD}, we assume that the sets $\ZadepsellNj$ have non-empty interior for all $j=0,\ldots,N_t$.

The only aspect to be explained is how to compute the snapshots in
Algorithm~\ref{Alg:MPCPOD}. We will explain this in the next
section. Let us mention that other possibilities are explored in
\cite{FFV19,GU14,Mec19,MV19}.

\section{Computational details and additional parallelization issues}
\label{sec:parallel}
In this section, we clarify some aspects of our simplification strategy to compute an an approximate solution to \eqref{eq:opt_cont_prob_in_u_bous}. We give details about the implementation of the algorithm, mentioning which parts can be performed in parallel. 

Looking at Algorithm~\ref{Alg:MPCPOD}, there are still some unclear aspects. At first, it is not clear how to compute the POD snapshots in Step~\ref{Alg:MPCPOD:step:compute_snapshots}. At iteration $j=0$, we initialize our model by solving the Boussinesq approximation \eqref{eq:Bous} with a given initial control $u_0$ for the time horizon $[0,T+N\Delta t]$ using the FE method for the spatial and the implicit Euler method for the temporal discretization. The intervall $[0,T+N\Delta t]$ is discretized by $t_j=j\Delta t$, $j=0,\ldots,m$ with $m=N_t+N$. Exploiting the power of parallel computing we obtain, therefore, the solution sequence $\{(\textbf{v}_{j,h}^0,y_{j,h}^0)\}_{j=0}^m$. Then, we compute the sequence $\{p_{j,h}^0\}_{j=0}^m$ by solving in parallel
\begin{align*}
&\frac{1}{\Delta t}\,{\langle p_{j,h}^0-p_{j+1,h}^0,\varphi_h\rangle}_H+a(t_j;\varphi_h,p_{j,h}^0)  +\frac{\sigma}{\varepsilon^2}\,{\langle \mathcal{G}^{0,\varepsilon}(t_j;y_{j,h}^0),\varphi_h\rangle}_H\\
&\qquad= \frac{\sigma}{\varepsilon^2}\,{\langle r^{0,\varepsilon}(t_j),\varphi_h\rangle}_H\quad\text{for all }\varphi_h\in V_h\text{ and }j=m-1,\ldots,1,\\
&\bigg\langle p_{m,h}^0 + \frac{\sigma}{\varepsilon^2}\,\mathcal{G}^{0,\varepsilon}(t_m;y_{m,h}^0),\varphi_h\bigg\rangle_H  = \frac{\sigma}{\varepsilon^2}\, {\langle r^{0,\varepsilon}(t_m),\varphi_h\rangle}_H\quad\text{for all }\varphi_h\in V_h
\end{align*}
(compare \eqref{eq:pdass_dual} and \eqref{pT}), where the terms $\mathcal{G}^{0,\varepsilon}$ and $r^{0,\varepsilon}$ involve the computation of initial active and inactive sets \eqref{ActiveSets}-\eqref{InactiveSets}, which can be also performed in parallel. The snapshots $\{y_{j,h}^0\}_{j=0}^m$ and $\{p_{j,h}^0\}_{j=0}^m$ are used to compute the POD basis in step~\ref{Alg:MPCPOD:step:compute_basis}. Note that the discrete POD eigenvalue problem has dimension $2m$, which may be large. We discuss later about this issue. For subsequent basis updates at an MPC step $j$, triggered by the a-posteriori error estimator, we will compute our snapshots with the same procedure, but using the control $\bar u^{\varepsilon,\ell}_N$ computed at the previous iteration $j-1$. Note that we get also an update for the velocity field $\textbf{v}$.

We use the velocity field $\textbf{v}$ (computed at 
step~\ref{Alg:MPCPOD:step:compute_snapshots}) as fixed one for the open-loop optimal control problem \eqref{eq:opt_control_problem_MPCPOD}, solved with POD-PDASS. In such a way, every time the basis and the velocity field are updated, the MPC can correct its approximation. We remark that the POD-PDASS matrix for problem \eqref{eq:opt_control_problem_MPCPOD} has the size $(2N\ell)\times (2N\ell)$ (and is thus small), therefore there is no need of parallelization for solving the linear system \eqref{DiscOpSystemPOD}. We anyway compute the active and inactive sets \eqref{ActiveSets}-\eqref{InactiveSets} in parallel, since they depend on $N_x$ and not on $\ell$. 

With the current information, one can already run Algorithm~\ref{Alg:MPCPOD}. In addition, we use a couple of improvements, which we describe further. A first disadvantage of the POD basis updates is that any time we need to compute the snapshots, we need to solve the Boussinesq equation for $t\in [t_j,T+N\Delta t]$. Although the measure of this time domain decreases for increasing $j$, computing the solution of the Boussinesq approximation \eqref{eq:Bous} for this large time horizon can be computationally demanding even using parallel computing. This computational effort can be reduced, considering to generate snapshots in a smaller time interval $[t_j,t_j+(N+M)\Delta t]$, where $M>0$ is chosen to be small. In such a way, we have an updated velocity field $\textbf{v}$ (and snapshots) for the subsequent $M$ iterations of the algorithm. If a POD basis update is performed in one of these steps, we have saved computational time. If not, we need to solve \eqref{eq:Bous} for the remaining time horizon up to the end (or for other $M$ steps, for example). Having only $2(N+M)$ snapshots is also beneficial for computing the POD basis, since we have to solve a smaller eigenvalue problem.

A second disadvantage of the POD basis update is that overwriting the past information while computing a new basis may not be always a good choice. The old snapshots, in fact, may still carry useful information for the future computations. Yet, adding only snapshots over snapshots will increase the number of POD basis needed to achieve the a-priori tolerance $\tau_1$. We employ, therefore, the snapshots selection strategy proposed in \cite{Mec19}. We report it in Algorithm~\ref{Alg:Snapshotselection}. This strategy consists in keeping only those old snapshots that represent dynamics close to the new ones. This procedure is motivated by the fact that the MPC algorithm gives a forecast of the future dynamics (new computed snapshots). Therefore, it makes sense to use them as starting point and add eventually previously computed information, if this describes similar trajectories. 
\begin{algorithm}
	\caption{Snapshots selection from \cite{Mec19}.\label{Alg:Snapshotselection}}
	\begin{algorithmic}[1]
		\REQUIRE Snapshots previously computed and stored in a list $\mathsf{L}$ and tolerances $0<\rho <\varrho \ll 1$.  
		\STATE Compute the new $2(N+M)$ snapshots and store them in a list $\mathsf{S}$;
		\FOR {$i \leq$ length($\mathsf{L}$)} 
		\FOR {$j \leq$ $2(N+M)$}
		\IF { \label{SnapSelIfcondition} $(1-\varrho)\|\mathsf{S}[j]\|\|\mathsf{L}[i]\| \leq |\langle \mathsf{S}[j],\mathsf{L}[i] \rangle | \leq (1-\rho) \|\mathsf{S}[j]\|\|\mathsf{L}[i]\| $}
		\STATE Add $\mathsf{L}[i]$ to $\mathsf{S}$;
		\STATE \textbf{break}
		\ENDIF
		\ENDFOR
		\ENDFOR
	\end{algorithmic}
\end{algorithm}
Furthermore, in step~\ref{Alg:MPCPOD:step:compute_feedback} of Algorithm \ref{Alg:MPCPOD}, one can solve in parallel one time step of the Boussinesq model \eqref{eq:Bous} applying the feedback law, instead of solving \eqref{StateEquation}. In such a way, we further improve the MPC approximation. This operation is not costly if one considers that it is done only for one time-step. At last, note that the computation of the a-posteriori error estimate \eqref{APostError} requires full-order solves of state and adjoint equations. This can also be carried out in parallel. To conclude, we resume every step of Algorithm~\ref{Alg:MPCPOD}, which is performed in parallel:
\begin{itemize}
	\item Computation of POD snapshots: this involves solving the Boussinesq approximation \eqref{eq:Bous} and the adjoint equation \eqref{eq:pdass_dual};
	\item Computation of the active and inactive sets \eqref{ActiveSets}-\eqref{InactiveSets} of the POD-PDASS;
	\item Computation of the initial guess for the next MPC iteration $j$: this involves solving one time step of \eqref{eq:Bous};
	\item The POD a-posteriori error estimator \eqref{APostError}.
\end{itemize}

\section{Numerical tests}
\label{sec:num_test}
The numerical tests were performed on the compute server of the Chair of Applied Mathematics at the University of Bayreuth, a computer with 32 physical cores (64 logical) Intel(R) Xeon(R) CPU E7-2830  @ 2.13GHz and a total RAM of 512GB. For the numerical realization of the Boussinesq approximation, we refer to \cite{And19}.
\subsection{Numerical verification of the turnpike property}
\label{sec:turnpike}
\begin{figure}
	\centering
	\subfigure[$\|y(t)\|_H$, with $y$ solution to \eqref{eq:opt_control_problem_MPC} for different horizons $N$.\label{Fig:open_loop_diff_horizon}]{\includegraphics[width= 60mm]{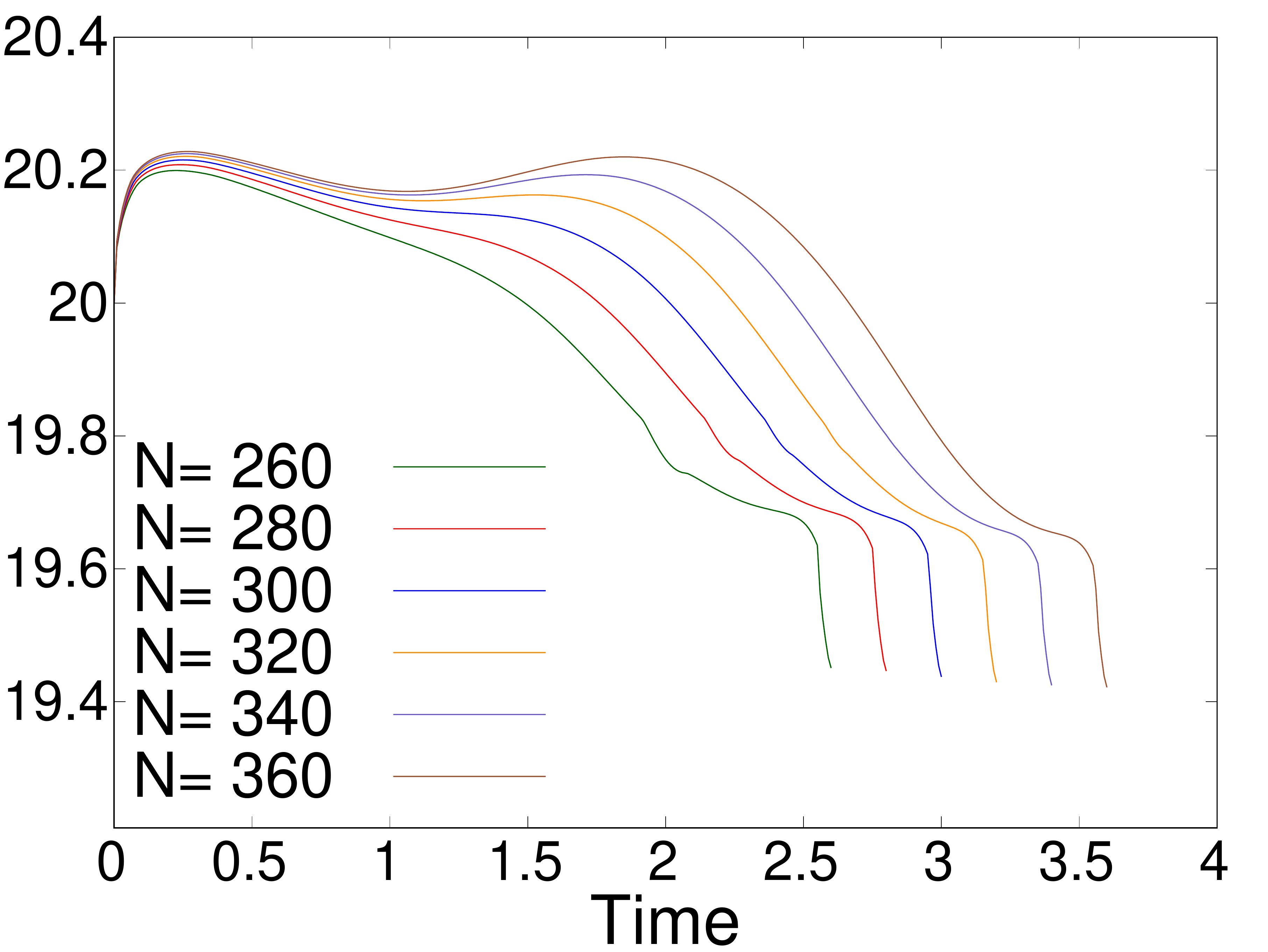}}
	\hspace{1em}
	\subfigure[$\|y(t)\|_H$, with $y$ solution to \eqref{eq:opt_control_problem_MPC} for different initial guesses and fixed horizon $N=360$. \label{Fig:open_loop_diff_init_guesses}]{\includegraphics[width=60mm]{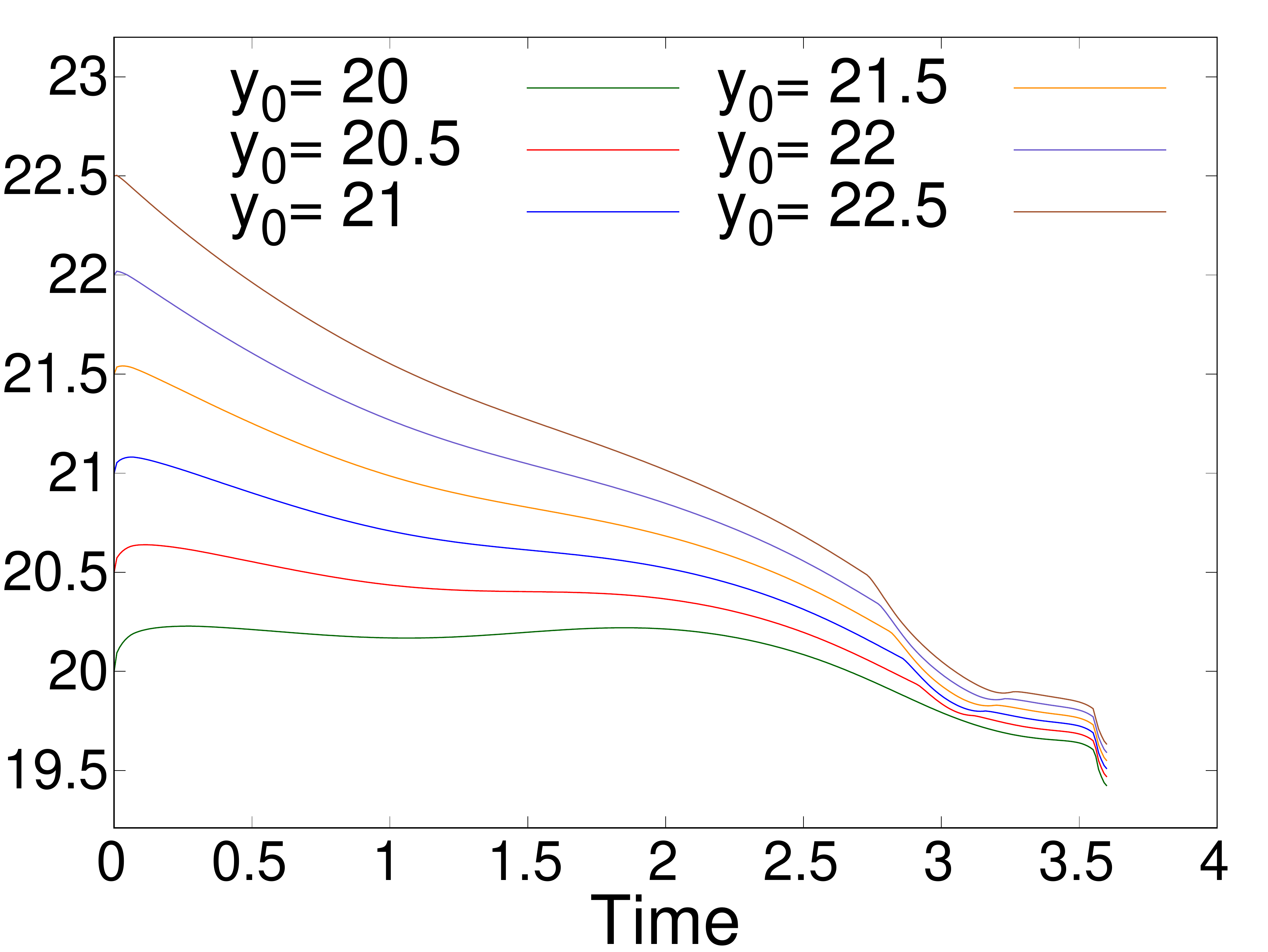}}
	\caption{Turnpike test results.}
\end{figure}

In this subsection we check numerically that the turnpike property holds for \eqref{eq:opt_control_problem_MPC}, because this is the most 
important ingredient for proving that MPC provides nearly optimal 
solutions, see \cite{GP20}. Before stating it in a mathematical formalism, we need to introduce the concept of optimal operation for an infinite horizon optimal control problem. Consider the discrete-time time-varying control system \eqref{eq:discrete_dynamical_system}. In this subsection, we denote with $y_u(\cdot\,;t_j,\tilde y_\circ)$ the solution of \eqref{eq:discrete_dynamical_system} for a given control $u \in \UNj$ starting in $\tilde y_\circ\in H$ at time $t_j$. Considering the stage cost function $l:\mathbb{R}\times H\to \mathbb{R}$, $l(u,w) = (|u|^2+\sigma\,\|w\|_H^2)/2$, the cost functional $J_N^{j}$ can be written as
\[
J_N^{j}(u,w) = \sum_{i=0}^N l(u(t_{j+i}),w(t_{j+i})).
\]
We recall that the finite horizon optimal control problem \eqref{eq:opt_control_problem_MPC} is 
\[
\min_{(u,w)\in \ZadepsNj} J_N^j(u,w)
\]
with (locally) optimal solution $(\bar u^\varepsilon_N,\bar w^\varepsilon_N)\in \ZadepsNj$. By $\Zadepsinfj$, $J_\infty^j$, $\bar u^\varepsilon_\infty$ and $\bar w^\varepsilon_\infty$ we denote the extension of the above quantities to infinite time horizon problem. Note that $\ZadepsNj$ changes if the starting value $\tilde y_\circ\in H$ of the PDE is different. Therefore, in what follows we indicate this dependency with $\ZadepsNj(\tilde y_\circ)$.
\begin{definition}
	\label{def:optimal_operation}
	Let $\tilde y_\circ\in H$ be an arbitrary feasible starting value at time $t_j$. Consider a control sequence $(\bar u,\bar w)\in \Zadepsinfj(\tilde y_\circ)$ with corresponding trajectory $\bar y$. We say that the system \eqref{eq:discrete_dynamical_system} is optimally operated at $(\bar y,\bar u,\bar w)$ if
	\[
	\liminf_{I\to\infty} \sum_{i=0}^{I-1} l(u(t_{j+i}),w(t_{j+i}))-l(\bar u(t_{j+i}),\bar w(t_{j+i}))\geq 0
	\]
	for all feasible $\hat y\in H$ and all $(u,w)\in \Zadepsinfj(\hat y)$. The trajectory $\bar y$ is called optimal operation trajectory.
\end{definition}
We remark that there may exist more than one solution to $\eqref{eq:discrete_dynamical_system}$ that satisfies Definition~\ref{def:optimal_operation}. Moreover, it is well-known that, in general, it is not easy to compute an optimal operation trajectory explicitly. Still, one can expect that the MPC trajectory converges to this trajectory for sufficiently large MPC horizon $N$. To ensure this behavior one has to assume structural conditions such as the turnpike property. This property demands that the solution to the finite (and infinite) time horizon optimal control problem is close for most of the time to an optimal operation trajectory $\bar y$, according to the following definition.
\begin{figure}
	\centering
	\includegraphics[]{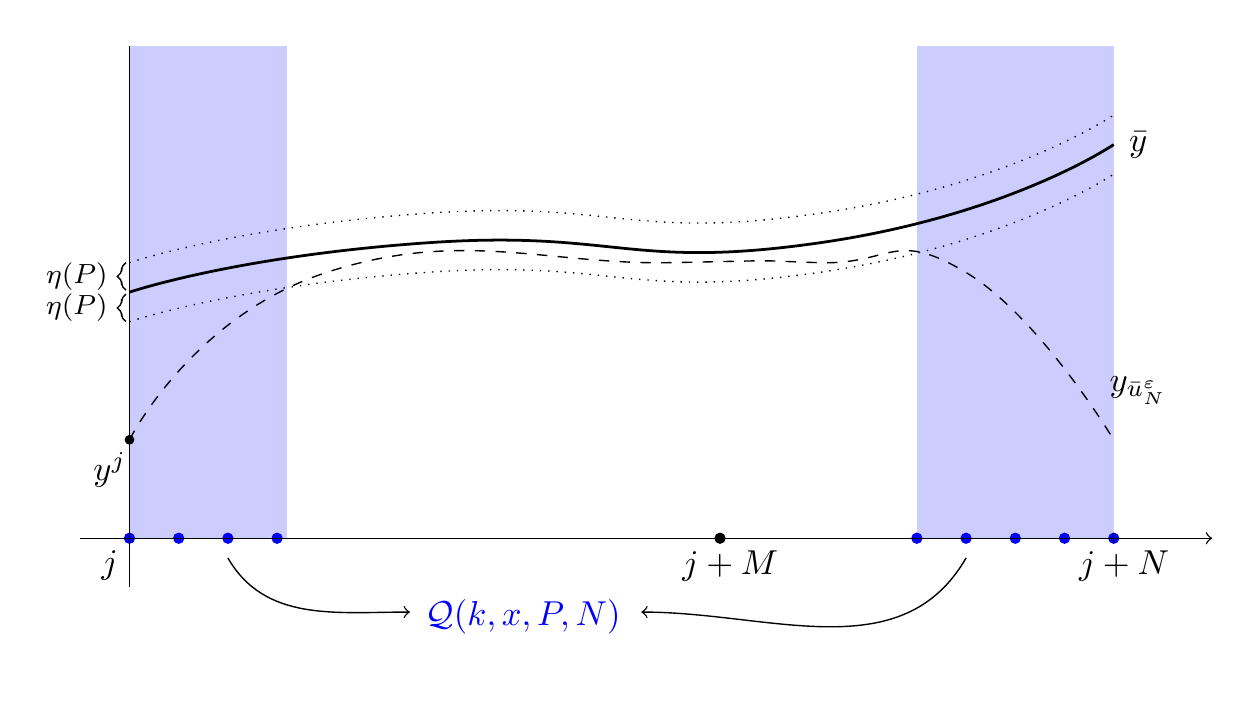}
	\caption{Time-varying turnpike property for the finite horizon optimal control problem. \label{Fig:turnpike}}
\end{figure}
\begin{definition}
	\label{def:Turnpike}
	Consider $(\bar y,\bar u,\bar w)$ at which system \eqref{eq:discrete_dynamical_system} is optimally operated. The finite horizon optimal control problem \eqref{eq:opt_control_problem_MPC} has the turnpike property at $(\bar y, \bar u, \bar w)$ if there exist $\eta\in \mathcal L\dagfootnote{ $\mathcal L= \{\eta\,:\,\mathbb{R}^+_0\to \mathbb{R}^+_0 \big| \eta \text{ is continuous and strictly decreasing with } \lim_{s\to\infty} \eta(s) = 0\}.$}$ such that for each $j\in\mathbb{N}_0$, each feasible initial guess $y_j$, each optimal trajectory $y_{\bar u^\varepsilon_N}(\cdot\,; t_j, y_j)$ and all $N,P\in\mathbb{N}$ there is a set $\mathcal{Q}(t_j,y_j,P,N)\subseteq\{0,\ldots,N\}$ with at most $P$ elements and
	\[
	\begin{aligned}
	& \|y_{\bar u^\varepsilon_N}(t_{j+M}; t_j, y^j)-\bar y(t_{j+M})\|_H+|\bar u^\varepsilon_N(t_{j+M})-\bar u(t_{j+M})| \\ & \quad + \|\bar w^\varepsilon_N(t_{j+M})-\bar w(t_{j+M})\|_H\leq \eta(P)
	\end{aligned}
	\]
	for all $M\in\{0,\ldots,N\}\setminus\mathcal{Q}(t_j,y_j,P,N)$.
\end{definition}
The definition of the turnpike property seems technically elaborated but it is easy to illustrate graphically; cf. Figure~\ref{Fig:turnpike}. This property demands that for each optimal open-loop trajectory there can only be a finite number $P$ of time steps, $P$ independent from $N$, at which the trajectory is far from the optimal operation trajectory. Moreover, it permits to define a bound on the distance for the remaining close points. Contrary to the computation of the optimal operation trajectory, it is easy to find numerical evidence that a given optimal control problem has the turnpike property \cite{GP19}.
\begin{figure}
	\centering
	\includegraphics[height=42mm]{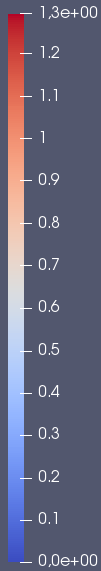}
	\includegraphics[height=42mm]{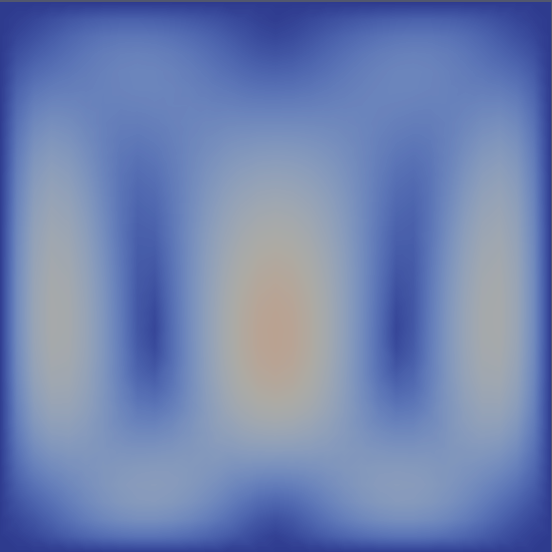}
	\includegraphics[height=42mm]{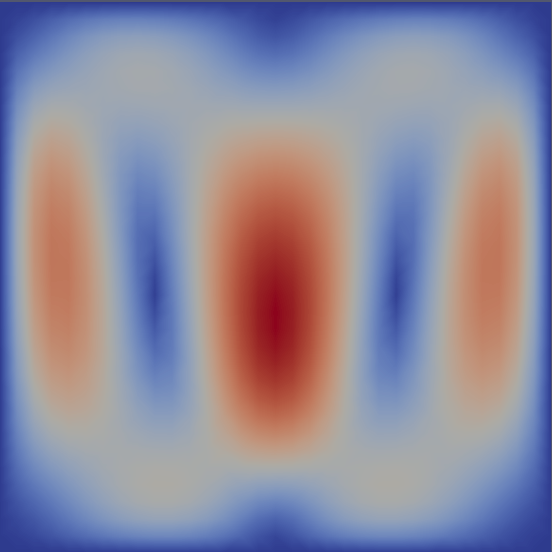}
	\includegraphics[height=42mm]{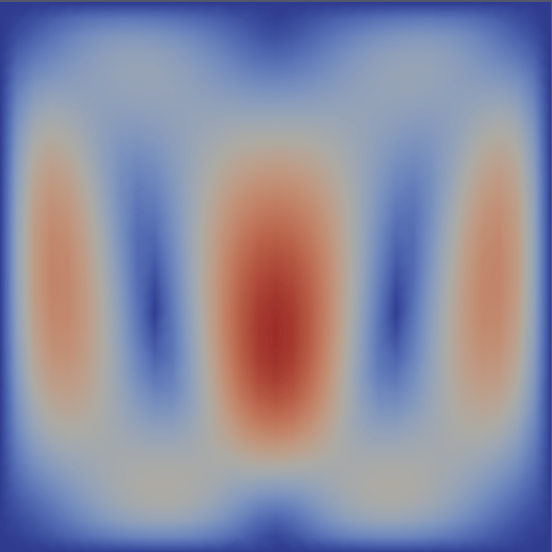}
	\caption{ Advection field \textbf{v} for the turnpike test at time instances $t= 1.0, 2.0, 3.0$. \label{Fig:advection_field_turnpike}}
\end{figure}
For this test, we choose $\Omega= (0,1)^2$ discretized with a structured quadrilateral mesh of $N_\bx= 1769$ nodes, $\Gamma_c= (0,1)\times \{0\}$, $\rho= C_p= \alpha= 1.0$, $\nu= 0.1$, $\kappa= 0.025$ and $\textbf{g} = (0,-9.8)$. Moreover, we consider $\gamma= 0.1$ and $\gamma_c= 10^5$ for the boundary conditions, so that the action of the control on the boundary can be seen (numerically) as the one of a Dirichlet boundary condition. We solve at first \eqref{eq:Bous} with fixed control $u(t)=22.5$, outside temperature $y_{\text{out}}= \max(13,16-t)$ and initial guess $(\textbf{v}_0,y_0)(t,\bx)=(0,20)$ for all $\bx\in\Omega$ and $t\in[0,T]$, in order to generate the advection field $\textbf{v}$ to use in \eqref{eq:opt_control_problem_MPC}; cf. Figure~\ref{Fig:advection_field_turnpike}. Then, we solve \eqref{eq:opt_control_problem_MPC} with $\ua(t) = 0$, $\ub(t) = 10^6$, $\wb= -\wa= 10^9$, $\varepsilon=0.025$, $\sigma=1.0$, $\ya(t) = 17.5+\min(t,2)$ and $\yb(t) = 23$ for all $t\in[0,T]$. We use a time step $\Delta t= 0.01$ and compute the open-loop problem solution for $j=0$ (i.e. $t_0=0$) for different horizons $N$, thus $T= N \Delta t$ in each test. As shown in Figure~\ref{Fig:open_loop_diff_horizon}, the $L^2$ norm of the solution to \eqref{eq:opt_control_problem_MPC} is the same for a larger number of time steps as $N$ increases. Moreover, the solutions show a similar behavior. 
This provides clear evidence for the turnpike behavior. This evidence is further strengthened by the results from Figure~\ref{Fig:open_loop_diff_init_guesses}, where the open-loop trajectories starting from different initial guess are converging to the same region as the time passes. For more information 
about how to check the conditions for near optimal performance of MPC 
numerically, we refer to \cite{GP19}.
\subsection{MPC-POD algorithm for the optimal control problem \eqref{eq:opt_cont_prob}}
\begin{figure}
	\centering
	\subfigure[Minimum (blue), maximum (red) and average (purple) temperature in the room and state constraints $\ya$ and $\yb$ (black). \label{Fig:MPCtest}]{\includegraphics[width= 0.47\textwidth]{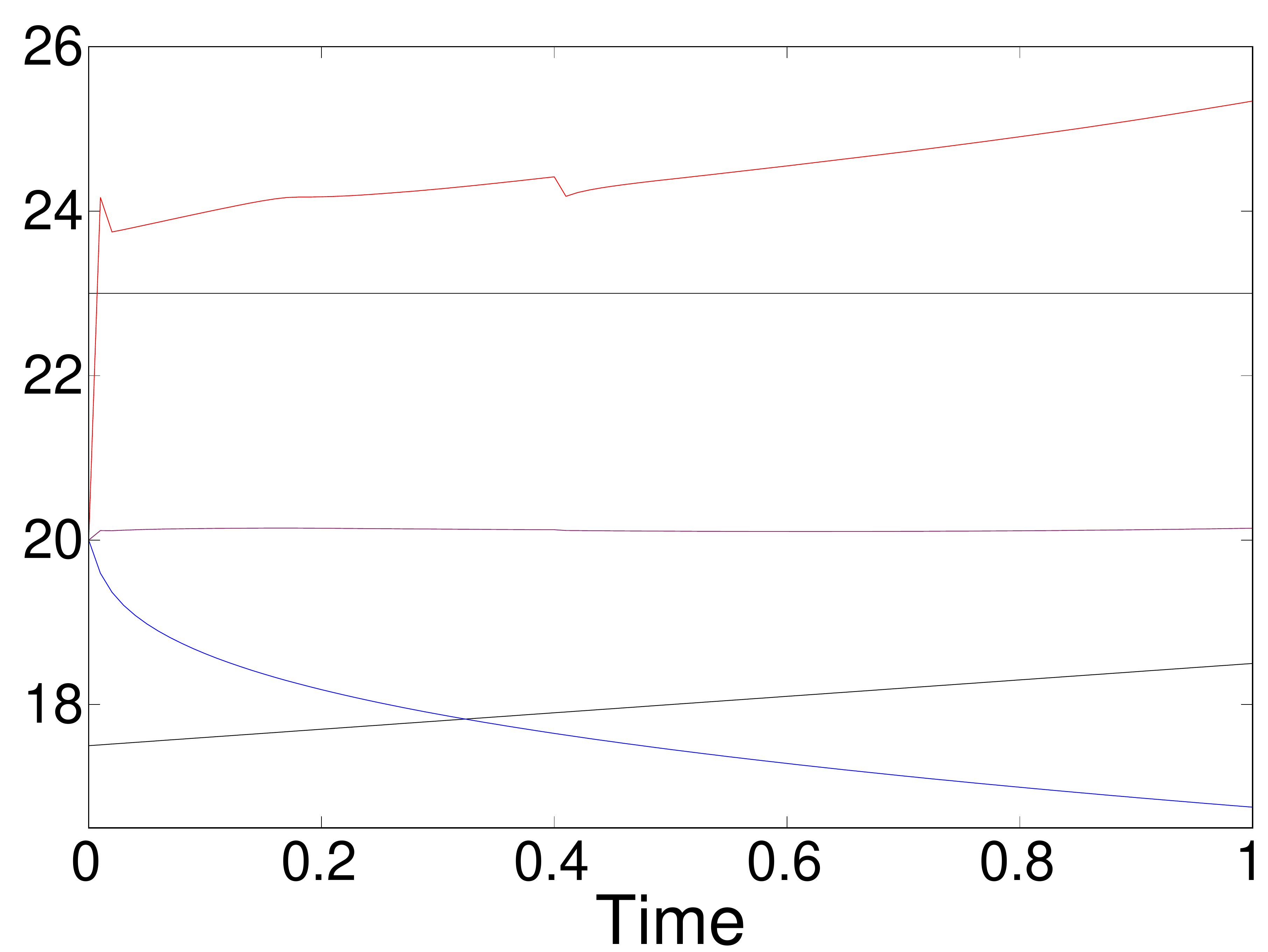}}
	\hspace{2pt}
	\subfigure[Number of active points with respect to the state constraints.\label{Fig:ActivePoints}]{\includegraphics[width= 0.47\textwidth]{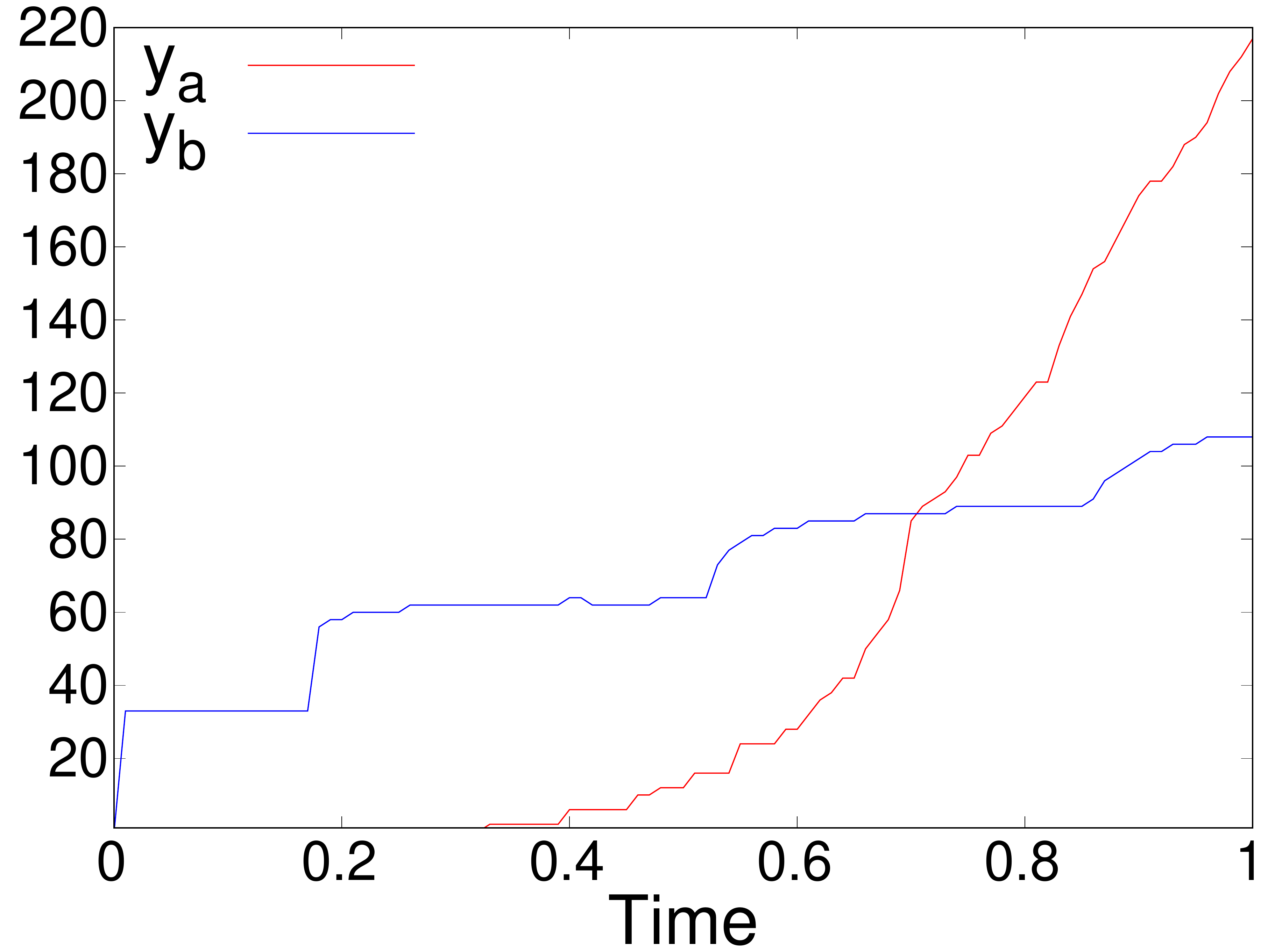}}
	\caption{Results for Algorithm~\ref{Alg:MPCPOD} with $\gamma= 0.1$ and $\varepsilon= 0.025$.}
\end{figure}
In this section we test Algorithm~\ref{Alg:MPCPOD} with the same data of Section~\ref{sec:turnpike}, choosing an MPC horizon $N=300$. This horizon is a good comprise between computational time and turnpike property; cf. Figure~\ref{Fig:open_loop_diff_horizon}. We fix also $M=60$, the a-priori tolerance $\tau_1= 10^{-8}\sum_{i=1}^\ell \lambda_i$ and the a posteriori error estimator tolerance $\tau_2= 3.5$. Although $\tau_2$ seems large, one has to consider that the control will have a value close to the state constraints range, roughly between 20 and 30, due to the fact that the control boundary condition is (numerically) close to a Dirichlet one. Moreover, the MPC horizon $N=300$ corresponds to a time horizon for the open-loop problem equal to three. Therefore, one can expect an average $L^2$ norm for the control of $43$ ($=25\sqrt{3}$). Hence, fixing $\tau_2= 3.5$ equals asking a relative error approximately of $8\%$. We expect, moreover, an overestimation factor for the a-posteriori error estimator, since $(\sigma^\mathsf{ap})^{-1}\|\mathcal{T}^\star \zeta\|_\U$ in \eqref{APostError} gives an upper-bound also for the error in reconstructing the artificial control, which might attain significantly greater values according to the chosen $\varepsilon$. In addition, a too small $\tau_2$ will trigger the basis updates too often and thus worsen the performance of the algorithm. 

Figure~\ref{Fig:MPCtest} shows the minimum, maximum and average value of the optimal MPC solution together with the desired bounds. One can notice that the average temperature is kept inside the bounds, where instead the minimum and the maximum one exceed them. This is typical of the $L^2$ regularization effect of the virtual control approach and can be mitigated taking a small $\varepsilon$. This effect is accentuated also by the POD method, which provides a sub-optimal solution as well-known from the literature \cite{Rav2000,TV09}. As a result, the average number of active points at each time step amounts to $3.5\%$ of the total number of grid points. This can be seen in Figure~\ref{Fig:ActivePoints}. Note that the number of points below the lower bound $\ya$ increases as the time passes. The reason is that the upper bound $\yb$ prevents the control at the boundary from growing significantly to counteract the decreasing outside temperature. This together with the sub-optimality of the POD solution will prevent the control from reacting faster and stronger to this external inputs. A possible solution is to consider local constraints in the domain $\Omega$ \cite{Pir20} and/or local POD basis functions \cite{AH2015,AHKO2012}. Nevertheless, the fact that the average temperature is inside the bounds and only a small portion of discretization points are active make Algorithm~\ref{Alg:MPCPOD} a viable approach in HVAC of residual buildings. 

\begin{figure}
	\centering
	\subfigure[Minimum (blue), maximum (red) and average (purple) temperature in the room and state constraints $\ya$ and $\yb$ (black). \label{Fig:MPCtest_smaller}]{\includegraphics[width= 0.47\textwidth]{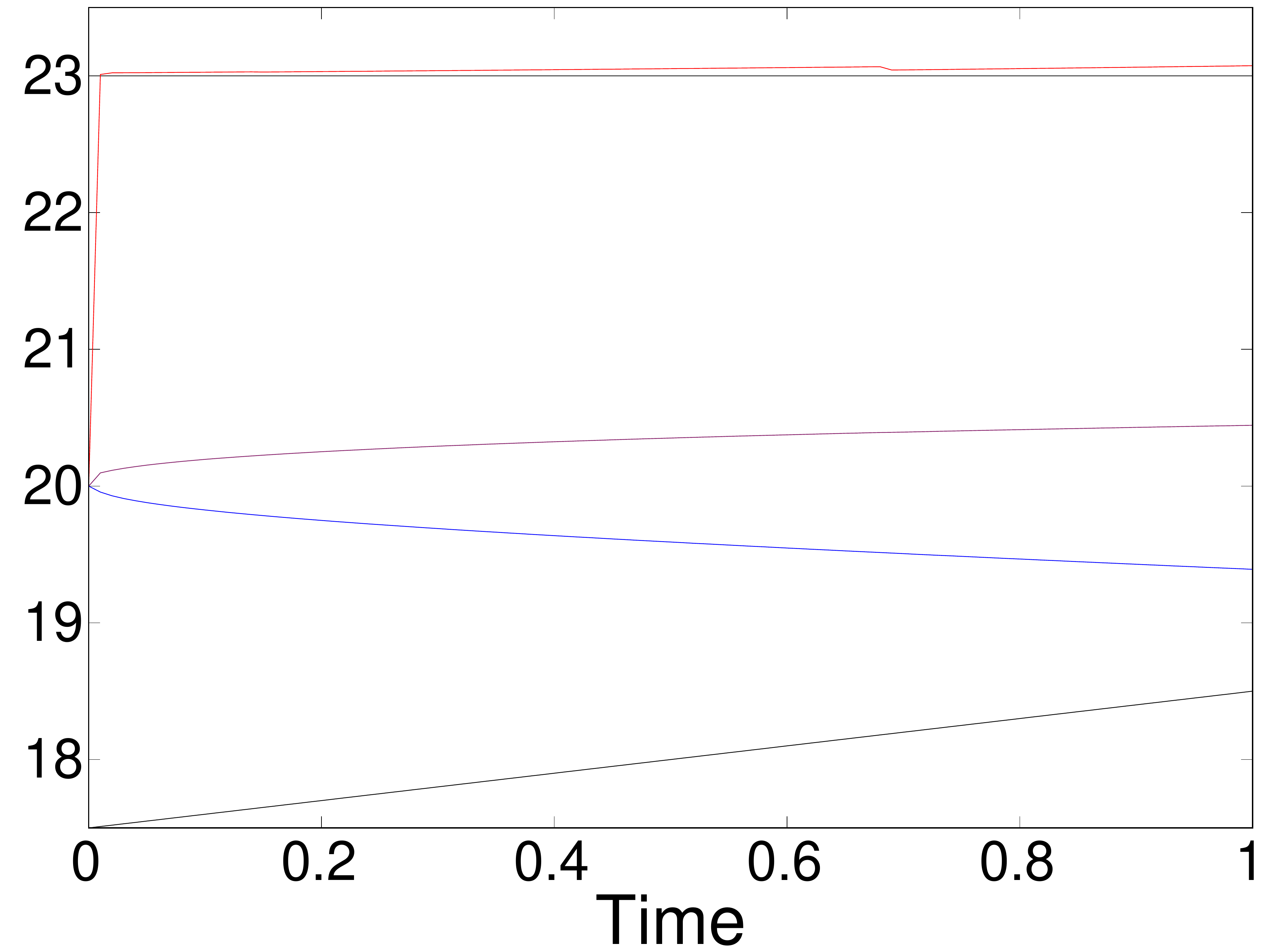}}
	\hspace{2pt}
	\subfigure[Number of active points with respect to the state constraints.\label{Fig:ActivePoints_smaller}]{\includegraphics[width= 0.47\textwidth]{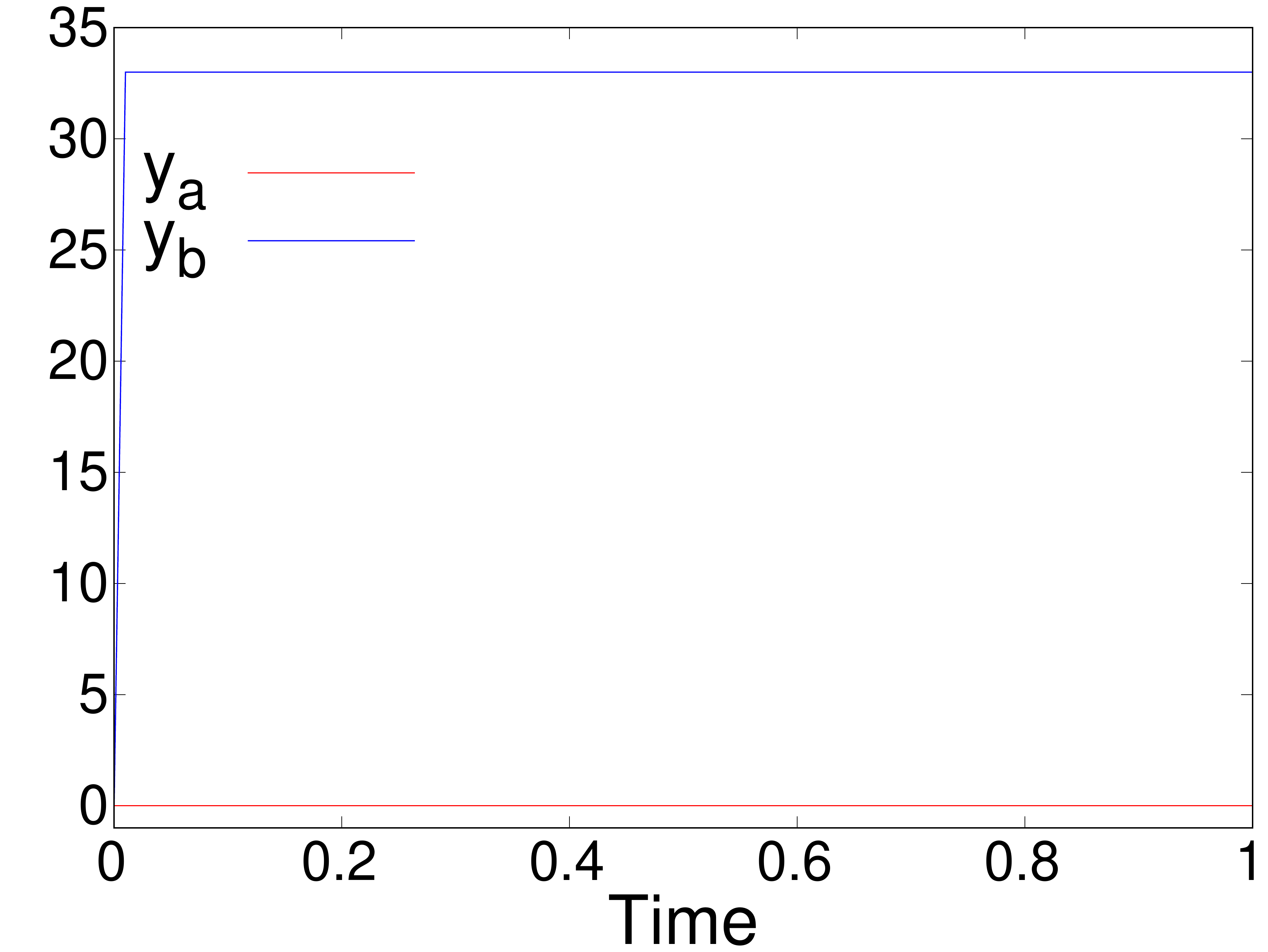}}
	\caption{Results for Algorithm~\ref{Alg:MPCPOD} with $\gamma= 0.01$ and $\varepsilon= 0.0075$. \label{Fig:smaller}}
\end{figure}

For smaller $\gamma=0.01$ and $\varepsilon=0.0075$, we can significantly improve the results in terms of violation of the state constraints;
cf. Figure~\ref{Fig:smaller}. In this test, in fact, the influence of the
outside temperature is lower and the method is capable to reduce
significantly the number of active points; see
Figure~\ref{Fig:ActivePoints_smaller}. In particular, the minimum
temperature in the room never lies below the given lower bound $\ya$
(Figure~\ref{Fig:MPCtest_smaller}). The maximum temperature, instead,
is slightly above the upper bound $\yb$. This can be justified from
the chosen value of the relaxation parameter $\varepsilon$ and also
from the fact that we are using the POD method, thus the reduced-order
model approximation introduce a further discretization error, which
may not be canceled even when passing to the limit $\varepsilon\to 0$. This
violation is anyway negligible for HVAC of residual buildings, since
the maximum value of the temperature attained all over the time-space
domain is $23.07$. Furthermore, we already discuss the conflict (and a
possible solution) between the state constraints and the control. This
can be clearly seen here, since the 33 active points are exactly the
grid points belonging to $\Gamma_c$.

\section{Conclusion}
We presented a POD-based economic MPC algorithm to handle HVAC of residual buildings. The strategy consists in computing a sub-optimal solution to an optimal control problem subjected to the Boussinesq approximation of Navier-Stokes and bilateral state and control constraints. The goal was to have a good compromise between robustness and computational cost. The first is achieved via the MPC and a goal-oriented a-posteriori error estimate, while the second is obtained through POD and parallel computing. The numerical results confirm the expectations. While some limitations of the strategy show up, particularly the emergence of several discretization points over time in which the state constraints are violated, the resulting control strategy can be regarded accurate enough for typical HVAC applications. The peformance can be further improved, for example, by considering localized MOR and localized constraints. This will be the focus of future work.

	\section*{Acknowledgments}
	The authors gratefully acknowledge support of DFG grants 274852524, 274852737 and 274853298. 
	This work was performed under the auspices of the U.S. Department of
	Energy by Lawrence Livermore National Laboratory under Contract DE-AC52-07NA27344
	(LLNL-JRNL-818244). This document was prepared as an account of work
	sponsored by an agency of the United States government. Neither the
	United States government nor Lawrence Livermore National Security, LLC,
	nor any of their employees makes any warranty, expressed or implied, or
	assumes any legal liability or responsibility for the accuracy,
	completeness, or usefulness of any information, apparatus, product, or
	process disclosed, or represents that its use would not infringe
	privately owned rights. Reference herein to any specific commercial
	product, process, or service by trade name, trademark, manufacturer, or
	otherwise does not necessarily constitute or imply its endorsement,
	recommendation, or favoring by the United States government or Lawrence
	Livermore National Security, LLC. The views and opinions of authors
	expressed herein do not necessarily state or reflect those of the United
	States government or Lawrence Livermore National Security, LLC, and
	shall not be used for advertising or product endorsement purposes.

	{\small
		\bibliographystyle{abbrv}
		\bibliography{bibliography}
	}
\end{document}